\documentclass[11pt]{article}
\openup 5pt\usepackage{mathrsfs}
\usepackage{latexsym,amsfonts,amssymb,amsthm}
\title{{\bf
Classification of quasifinite representations with nonzero central
charges for type $A_1$ EALA with coordinates in quantum
torus}\thanks{\ Supported by the National Science Foundation of
China (No. 10671160, 10471091), the China Postdoctoral Science
Foundation (No. 20060390693), and ``One Hundred Talents Program''
from University of Science and Technology of China. } }
\author{   Weiqiang\ Lin$^{1,2}$
\quad  Yucai\ Su$^2$ \\
 {\small 1. Department of Mathematics, Zhangzhou
Teachers College,}\\
{\small Zhangzhou\ 363000, Fujian, China}\\
{\small 2. Department of Mathematics, University of Science and Technology of China,}\\
{\small Hefei\ 230026, Anhui, China} \\ {\small E-mail:
linwq83@yahoo.com.cn, ycsu@ustc.edu.cn} }
\date{}
\def\la{\langle}
\def\ra{\rangle}
\begin{document}
\maketitle\vskip-20pt

 {\small {\bf Abstract:} In this paper, we first construct a Lie
 algebra $L$ from rank 3 quantum torus, and show that it is isomorphic
 to the core of EALAs of type $A_1$ with coordinates in rank
2 quantum torus. Then we construct two classes of irreducible ${\bf
Z}$-graded highest weight representations, and give
 the necessary and sufficient conditions for these representations to be quasifinite.
 Next, we prove that they
 exhaust all the generalized highest weight irreducible ${\bf
 Z}$-graded quasifinite representations. As a consequence, we determine all the irreducible ${\bf
 Z}$-graded quasifinite representations
 with nonzero central charges. Finally,
 we construct two classes of highest weight ${\bf
 Z}^2$-graded quasifinite representations by using these ${\bf Z}$-graded modules.

 {\bf Keyword:} core of EALAs, graded representations, quasifinite representations,
 highest weight representations, quantum
 torus.\vskip7pt
}

\centerline{\large\bf \S 1\ \  Introduction}

Extended affine Lie algebras (EALAs) are higher dimensional
generalizations of affine Kac-Moody Lie algebras introduced in [1]
(under the name of irreducible quasi-simple Lie algebras). They can
be roughly described as complex Lie algebras which have a
nondegenerate invariant form, a self-centralizing finite-dimensional
ad-diagonalizable Abelian subalgebra (i.e., a Cartan subalgebra), a
discrete irreducible root system, and ad-nilpotency of nonisotropic
root spaces (see [2--4]). Toroidal Lie algebras, which are universal
central extensions of $\dot{\mathfrak{g}}\otimes{\bf C}[t_1^{\pm
1},\cdots,t_n^{\pm 1}]$ ($\dot{\mathfrak{g}}$ is a
finite-dimensional simple Lie algebra), are prime examples of EALAs
studied in [5--11], among others. There are many EALAs which allow
not only Laurent polynomial algebra ${\bf C}[t_1^{\pm 1}, \cdots,
t_n^{\pm 1}]$ as coordinate algebra but also quantum tori, Jordan
tori and the octonians tori as coordinated algebras depending on
type of the Lie algebra (see [2, 3, 12--14]). The structure theory
of the EALAs of type $A_{d-1}$ is tied up with Lie algebra
$gl_d({\bf C})\otimes{\bf C}_Q$ where ${\bf C}_Q$ is the quantum
torus. Quantum torus defined in [15] are noncommutative analogue of
Laurent polynomial algebras. The universal center extension of the
derivation Lie algebra of rank 2 quantum torus is known as the
$q$-analog Virasoro-like algebra (see [16]). Representations for Lie
algebras coordinated  by certain quantum tori have been studied by
many people (see [17--22] and the references therein). The structure
and representations of the $q$-analog Virasoro-like algebra are
studied in many papers (see [23--27]).
%There are many relations between the representations of these two algebras.
In this paper, we first construct a Lie algebra $L$ from rank 3
quantum torus, which contains the $q$-analog Virasoro-like algebra
as its Lie subalgebra, and show that it is isomorphic to the core of
EALAs of type $A_1$ with coordinates in rank 2 quantum torus. Then
we study quasifinite representation of $L$.

When we study quasifinite representations of a Lie algebra of this
kind, as pointed out by Kac and Radul in [28], we encounter the
difficulty that though it is {\bf Z}-graded, the graded subspaces
are still infinite dimensional, thus the study of quasifinite
modules is a nontrivial problem.

Now we explain this paper in detail. In Section 2, we first recall
some concepts about the quantum torus and EALAs of type $A_1$. Next,
we construct a Lie algebra $L$ from a special class of rank 3
quantum, and show that $L$ is isomorphic to the core of EALAs of
type $A_1$ with coordinates in rank 2 quantum torus. Then, we prove
some basic propositions and reduce the classification of irreducible
${\bf Z}$-graded representations of $L$ to that of the generalized
highest weight representations and the uniformly bounded
representations. In Section 3, we construct two class of irreducible
${\bf Z}$-graded highest weight representations of $L$, and give
 the necessary and sufficient conditions for these representations to be quasifinite.
 In Section 4, we prove
 that the generalized highest weight irreducible ${\bf Z}$-graded
quasifinite representations of  $L$  must be
 the highest weight representations, and thus the representations constructed in
 Section 3 exhaust all the generalized highest weight quasifinite representations. As a
 consequence, we complete the classification of irreducible ${\bf
 Z}$-graded quasifinite representations of $L$ with nonzero central charges,
 see Theorem 4.4 (the Main Theorem). In Section 5, we
 construct two classes of highest weight ${\bf Z}^2$-graded quasifinite representations.
\vskip7pt

\centerline{\large\bf \S 2\ \  Basics}

 Throughout this paper we use ${\bf
C}, {\bf Z}, {\bf Z}_+, {\bf N}$ to denote the sets of complex
numbers, integers, nonnegative integers, positive integers
respectively. And we use ${\bf C}^*, {\bf Z}^{2*}$ to denote the set
of nonzero complex numbers and ${\bf Z}^2\setminus\{(0,0)\}$
respectively. All spaces considered in this paper are over ${\bf
C}$. As usual, if $u_1,u_2,\cdots,u_k$ are elements on some vector
space, we use $\la u_1,\cdots,u_k\ra $ to denote their linear span
over ${\bf C}$. Let $q$ be a nonzero complex number. We shall fix a
generic $q$ throughout this paper.

\vspace{2mm} Now we recall the concept of quantum torus from [15].
Let $\nu$ be a positive integer and $Q=(q_{ij})$ be a $\nu\times\nu$
matrix, where
$$
q_{ij}\in{\bf C}^*,\; q_{ii}=1,\;q_{ij}=q_{ji}^{-1},\quad
\mbox{for}\; 0\leq i,j\leq \nu-1.
$$
A quantum torus associated to $Q$ is the unital associative algebra
${\bf C}_Q[t_0^{\pm 1},\cdots,t_{\nu-1}^{\pm 1}]$ (or, simply ${\bf
C}_Q$) with generators $t_0^{\pm 1},\cdots,t_{\nu-1}^{\pm 1}$ and
relations
$$
t_it_i^{-1}=t_i^{-1}t_i=1\;\mbox{and}\;t_it_j=q_{ij}t_jt_i,\quad\forall\
0\leq i,j\leq\nu-1.
$$
Write $t^{\bf m}=t_0^{m_0}t_1^{m_1}\cdots t_{\nu-1}^{m_{\nu-1}}$ for
${\bf m}=(m_0,m_1,\cdots,m_{\nu-1})$. Then
$$
t^{\bf m}\cdot t^{\bf n}=\Big(\,\mbox{$\prod\limits_{0\leq j\leq
i\leq{\nu-1}}$}q_{ij}^{m_in_j}\Big)t^{{\bf m}+{\bf n}},
$$
where ${\bf m, n}\in{\bf Z}^{\nu}$. If $Q=\Big(\begin{array}{cc} 1 &
q^{-1} \\ q & 1 \end{array}\Big)$, we will simply denote ${\bf C}_Q$
by ${\bf C}_q$.

\vspace{2mm}
Next we recall the construction of EALAs of type $A_1$
with coordinates in ${\bf C}_{q^2}$. Let $E_{ij}$ be the $2\times 2$
matrix which is $1$ in the $(i,j)$-entry and 0 everywhere else. The
Lie algebra $\widetilde{\tau}=gl_2({\bf C}_{q^2})$ is defined by
$$
[E_{ij}(t^{\bf m}),E_{kl}(t^{\bf
n})]_0=\delta_{j,k}q^{2m_2n_1}E_{il}(t^{{\bf m}+{\bf
n}})-\delta_{l,i}q^{2n_2m_1}E_{kj}(t^{{\bf m}+{\bf n}}),
$$
where $1\leq i,j,k,l\leq 2, {\bf m}=(m_1,m_2)$ and ${\bf
n}=(n_1,n_2)$ are in ${\bf Z}^2$. Thus the derived Lie subalgebra of
$\widetilde{\tau}$ is $\overline{\tau}=sl_2({\bf C}_{q^2})\oplus \la
I(t^{\bf m})\mid {\bf m}\in{\bf Z}^{2*}\ra $, where
$I=E_{11}+E_{22}$, since $q$ is generic. And the universal central
extension of $\overline{\tau}$ is $\tau=\overline{\tau}\oplus\la
K_1,K_2\ra $ with the following Lie bracket
$$
[X(t^{\bf m}),Y(t^{\bf n})]=[X(t^{\bf m}),Y(t^{\bf
n})]_0+\delta_{{\bf m}+{\bf n},0}q^{2m_2n_1}(X,Y)(m_1K_1+m_2K_2),
$$$$
\mbox{ }\; K_1,K_2 \;\mbox{are central},
$$
where $X(t^{\bf m}),Y(t^{\bf n})\in\overline{\tau}$ and $(X,Y)$ is
the trace of $XY$. The Lie algebra $\tau$ is the core of the EALAs
of type $A_1$ with coordinates in ${\bf C}_{q^2}$. If we add degree
derivations $d_1, d_2$ to $\tau$, then $\tau\oplus\la d_1,d_2\ra $
becomes an EALAs since $q$ is generic.

\vspace{2mm}
Now we construct our Lie algebra. Let
$$
Q=\left(\begin{array}{lll}
1 & -1 & 1 \\
-1 & 1 & q^{-1}\\
1 & q & 1
\end{array}\right).
$$
Let $J$ be the two-sided ideal of ${\bf C}_Q$ generated by
$t_0^2-1$. Define
$$
\widetilde{L}={\bf C}_Q/J=\la t_0^{i}t_1^{j}t_2^{k}\mid i\in{\bf
Z}_2,\ j,k\in{\bf Z}\ra ,
$$
be the quotient of ${\bf C}_Q$ by $J$ and identify $t_0$ with its
image in $\widetilde{L}$. Then the derived Lie subalgebra of
$\widetilde{L}$ is $\overline{L}=\la t_0^{\bar 0}t^{\bf m}\mid {\bf
m}\in{\bf Z}^{2*}\ra \oplus\la t_0^{\bar 1}t^{\bf m}\mid {\bf
m}\in{\bf Z}^2\ra $. Now we define a central extension of
$\overline{L}$, which will be denoted by $L=\overline{L}\oplus\la
c_1,c_2\ra $, with the following Lie bracket
$$
[t_0^it^{\bf m},t_0^jt^{\bf n}]=
((-1)^{m_1j}q^{m_2n_1}-(-1)^{in_1}q^{m_1n_2})t_0^{i+j}t^{{\bf
m}+{\bf n}}+ (-1)^{m_1j}q^{m_2n_1}\delta_{i+j,{\bar 0}}\delta_{{\bf
m}+{\bf n},0}(m_1c_1 +m_2c_2),
$$
$$
\mbox{ }\;\;c_1,c_2\;\; \mbox{are central,}
$$
where $i,j\in{\bf Z}_2$, ${\bf m}=(m_1,m_2)$ and ${\bf n}=(n_1,n_2)$
are in ${\bf Z}^2$. One can easily see that $\la t_0^{\bar 0}t^{\bf
m}\mid {\bf m}\in{\bf Z}^{2*}\ra \oplus\la c_1,c_2\ra $ is a Lie
subalgebra of $L$, which is isomorphic to the $q$-analog
Virasoro-like algebra.

 First we prove that the Lie algebra $L$ is in fact
isomorphic to the core of the EALAs of type $A_1$ with coordinates
in ${\bf C}_{q^2}$.\vspace{3mm}

 {\em {\bf Proposition 2.1 } The Lie algebra $L$
is isomorphic to $\tau$ and the isomorphism is given by the linear
extension of the following map $\varphi$:
\begin{eqnarray*}
t_0^{i}t_1^{2m_1+1}t_2^{m_2}&\mapsto&
(-1)^iq^{-m_2}E_{12}(t_1^{m_1}t_2^{m_2})+E_{21}(t_1^{m_1+1}t_2^{m_2}),
\\
t_0^{i}t_1^{2m_1}t_2^{m_2}&\mapsto&
(-1)^iE_{11}(t_1^{m_1}t_2^{m_2})+q^{-m_2}E_{22}(t_1^{m_1}t_2^{m_2})+\delta_{i,\bar
1}\delta_{m_1,0}\delta_{m_2,0}\frac{1}{2}K_1,
\\
c_1&\mapsto& K_1,\quad\ \quad\ \quad \ \quad c_2\ \ \mapsto \ \
2K_2,
\end{eqnarray*}
where $t_0^{i}t_1^{2m_1+1}t_2^{m_2},\ t_0^{i}t_1^{2m_1}t_2^{m_2}\in
L$.}\vspace{2mm}

{\bf Proof } We need to prove that $\varphi$ preserves Lie bracket.
First we have
$$
\begin{array}{ll}
[(-1)^iq^{-m_2}E_{12}(t_1^{m_1}t_2^{m_2})+E_{21}(t_1^{m_1+1}t_2^{m_2}),
(-1)^jq^{-n_2}E_{12}(t_1^{n_1}t_2^{n_2})+E_{21}(t_1^{n_1+1}t_2^{n_2})]
\\[7pt]
=
\Big((-1)^jq^{m_2(2n_1+1)}-(-1)^iq^{n_2(2m_1+1)}\Big)\Big((-1)^{i+j}E_{11}(t_1^{m_1+n_1+1}t_2^{m_2+n_2})
\\[7pt]
\phantom{=\big((-1)^jq^{m_2(2n_1+1)}-(-1)^iq^{n_2(2m_1+1)}\big)\Big(}
+q^{-m_2-n_2}E_{22}(t_1^{m_1+n_1+1}t_2^{m_2+n_2})\Big)\\[7pt]
\phantom{=}+\delta_{m_1+n_1+1,0}\delta_{m_2+n_2,0}(-1)^jq^{m_2(2n_1+1)}\Big((-1)^{i+j}(m_1K_1+m_2K_2)+(m_1+1)K_1+m_2K_2\Big)
\\[7pt]
=\Big((-1)^jq^{m_2(2n_1+1)}-(-1)^iq^{n_2(2m_1+1)}\Big)\Big((-1)^{i+j}E_{11}(t_1^{m_1+n_1+1}t_2^{m_2+n_2})
\\[7pt]\phantom{=\Big((-1)^jq^{m_2(2n_1+1)}-(-1)^iq^{n_2(2m_1+1)}\Big)\Big(}
+q^{-m_2-n_2}E_{22}(t_1^{m_1+n_1+1}t_2^{m_2+n_2})\Big)
\\[7pt]
\phantom{=}+\delta_{i+j,\bar
0}\delta_{m_1+n_1+1,0}\delta_{m_2+n_2,0}(-1)^jq^{m_2(2n_1+1)}((2m_1+1)K_1+2m_2K_2)
\\[7pt]\phantom{=}
 +\delta_{i+j,\bar
1}\delta_{m_1+n_1+1,0}\delta_{m_2+n_2,0}(-1)^jq^{m_2(2n_1+1)}K_1.
\end{array}
$$
On the other hand, we have
$$
\begin{array}{ll}
[t_0^it_1^{2m_1+1}t_2^{m_2},t_0^jt_1^{2n_1+1}t_2^{n_2}]
=\!\!\!\!&\Big((-1)^{j}q^{m_2(2n_1+1)}-(-1)^iq^{(2m_1+1)n_2}\Big)t_0^{i+j}
t_1^{2m_1+2n_1+2}t_2^{m_2+n_2} \\[7pt]
 &+\delta_{i+j,\bar
0}\delta_{2m_1+2n_1+2,0}\delta_{m_2+n_2,0}(-1)^jq^{m_2(2n_1+1)}((2m_1+1)c_1+m_2c_2).
\end{array}$$ Thus
$$\varphi([t_0^{i}t_1^{2m_1+1}t_2^{m_2}),t_0^{j}t_1^{2n_1+1}t_2^{n_2}])=[\varphi(t_0^{i}t_1^{2m_1+1}t_2^{m_2}),
\varphi(t_0^{j}t_1^{2n_1+1}t_2^{n_2})].$$  Similarly, we have
$$
\begin{array}{ll}
[\varphi(t_0^{i}t_1^{2m_1}t_2^{m_2}),\varphi(t_0^{j}t_1^{2n_1}t_2^{n_2})]
\\[7pt]
=[(-1)^iE_{11}(t_1^{m_1}t_2^{m_2})+q^{-m_2}E_{22}(t_1^{m_1}t_2^{m_2}),(-1)^jE_{11}(t_1^{n_1}t_2^{n_2})+q^{-n_2}E_{22}(t_1^{n_1}t_2^{n_2})
]
\\[7pt]
=(q^{2m_2n_1}-q^{2n_2m_1})\Big
((-1)^{i+j}E_{11}(t_1^{m_1+n_1}t_2^{m_2+n_2})+q^{-m_2-n_2}E_{22}(t_1^{m_1+n_1}t_2^{m_2+n_2})\Big)
\\[7pt]\phantom{=}
+\delta_{m_1+n_1,0}\delta_{m_2+n_2,0}\delta_{i+j,\bar
0}q^{2m_2n_1}(2m_1K_1+2m_2K_2), \end{array}$$ and
$$
\begin{array}{ll}
[t_0^{i}t_1^{2m_1}t_2^{m_2},t_0^{j}t_1^{2n_1}t_2^{n_2}]
=\!\!\!\!&(q^{2m_2n_1}-q^{2m_1n_2})t_0^{i+j}t_1^{2m_1+2n_1}t_2^{m_2+n_2}\\[7pt]
&+\delta_{i+j,\bar
0}\delta_{m_1+n_1,0}\delta_{m_2+n_2,0}q^{2m_2n_1}(2m_1c_1+m_2c_2).
\end{array}$$ Therefore
$$[\varphi(t_0^{i}t_1^{2m_1}t_2^{m_2}),\varphi(t_0^{j}t_1^{2n_1}t_2^{n_2})]=
\varphi([t_0^{i}t_1^{2m_1}t_2^{m_2},t_0^{j}t_1^{2n_1}t_2^{n_2}]).$$
Finally, we have
$$
\begin{array}{ll}
[\varphi(t_0^{i}t_1^{2m_1+1}t_2^{m_2}),\varphi(t_0^{j}t_1^{2n_1}t_2^{n_2})]
\\[7pt]
=[(-1)^iq^{-m_2}E_{12}(t_1^{m_1}t_2^{m_2})+E_{21}(t_1^{m_1+1}t_2^{m_2}),
(-1)^jE_{11}(t_1^{n_1}t_2^{n_2})+q^{-n_2}E_{22}(t_1^{n_1}t_2^{n_2})]
\\[7pt]
=\Big((-1)^jq^{2m_2n_1}-q^{n_2(2m_1+1)}\Big)\Big((-1)^{i+j}q^{-m_2-n_2}E_{12}(t_1^{m_1+n_1}t_2^{m_2+n_2})
 +E_{21}(t_1^{m_1+n_1+1}t_2^{m_2+n_2})\Big),
\end{array}$$ and
$$
[t_0^{i}t_1^{2m_1+1}t_2^{m_2},t_0^{j}t_1^{2n_1}t_2^{n_2}]
=((-1)^{j}q^{2m_2n_1}-q^{n_2(2m_1+1)})t_0^{i+j}t_1^{2m_1+2n_1+1}t_2^{m_2+n_2}.
$$
Thus
$$[\varphi(t_0^{i}t_1^{2m_1+1}t_2^{m_2}),\varphi(t_0^{j}t_1^{2n_1}t_2^{n_2})]=
\varphi([t_0^{i}t_1^{2m_1+1}t_2^{m_2},t_0^{j}t_1^{2n_1}t_2^{n_2}])
.$$ This completes the proof. \hfill$\Box$

\vspace{2mm} {{\bf Remark 2.2 } From the proof of above proposition,
one can easily see that $gl_2({\bf C}_{q^2})\cong\widetilde{L}$ and
$\overline{\tau}\cong\overline{L}$.}

\vspace{2mm} Next we will recall some concepts about the ${\bf
Z}$-graded $L$-modules. Fix a ${\bf Z}$-basis $${\bf
m}_1=(m_{11},m_{12}),\ {\bf m}_2=(m_{21},m_{22})\in{\bf Z}^2.$$ If
we define the degree of the elements in $\la t_0^it^{j{\bf m}_1
+k{\bf m}_2}\in L\mid i\in{\bf Z}_2,k\in{\bf Z}\ra $ to be $j$ and
the degree of the elements in $\la c_1,c_2\ra $ to be zero, then $L$
can be regarded as a ${\bf Z}$-graded Lie algebra:
$$
L_j=\la t_0^{i}t^{j{\bf m}_1+k{\bf m}_2}\in L\mid i\in{\bf
Z}_2,k\in{\bf Z}\ra \oplus\delta_{j,0}\la c_1,c_2\ra .
$$Set$$
L_+=\bigoplus\limits_{j\in{\bf N}}L_j,\quad
L_-=\bigoplus\limits_{-j\in{\bf N}}L_j.
$$
 Then $L=\oplus_{j\in{\bf
Z}}L_j$ and $L$ has the following triangular decomposition
$$
L=L_-\oplus L_0\oplus L_+.
$$\vskip2mm

{\bf Definition } \def\Z{{\bf Z}} For any $L$-module $V$, if
$V=\oplus_{m\in {\bf Z}} V_m$ with
$$
L_j\cdot V_m\subset V_{m+j},\;\forall\ j,m\in{\bf Z},
$$
 then $V$ is called a {\it {\bf Z}-graded
$L$-module} and $V_{m}$ is called a {\it homogeneous subspace of $V$
with degree $m\in{\bf Z}$.} The $L$-module $V$ is called
\begin{itemize}\parskip-1pt\item[(i)]
a {\em quasi-finite {\bf Z}-graded module} if ${\rm
dim}\,V_m<\infty,\,\forall\,m\in\Z$;
\item[(ii)] a {\em  uniformly bounded module}
if there exists some $N\in{\bf N}$ such that ${\rm
dim}\,V_m\le N,\,\forall\,m\in\Z$;
\item[(iii)] a {\em highest {\rm(resp.~}lowest$)$ weight module} if there exists a nonzero homogeneous
vector $v\in V_m$ such that $V$ is generated by $v$ and $L_+\cdot
v=0$ (resp.~$L_-\cdot v=0$);
\item[(iv)] a
{\em generalized highest weight module with highest degree $m$}
(see, e.g., [31]) if there exist a {\bf Z}-basis $B=\{{\bf b_1},{\bf
b_2}\}$ of ${\bf Z}^2$ and a nonzero vector $v\in V_{m}$ such that
$V$ is generated by $v$ and $t_0^it^{\bf m}\cdot v=0,\forall\ {\bf
m}\in{\bf Z}_{+}{\bf b_1}+{\bf Z}_{+}{\bf b_2},i\in{\bf Z}_2$;
\item[(v)] an {\em  irreducible {\bf Z}-graded module} if $V$ does not
have any nontrivial $\bf Z$-graded submodule (see, e.g., [29]).
\end{itemize}

We denote the set of quasi-finite irreducible ${\bf Z}$-graded
$L$-modules by ${\cal O}_{\bf Z}$. From the definition, one sees
that the generalized highest weight modules contain the highest
weight modules and the lowest weight modules as their special cases.
As the central elements $c_1,\ c_2$ of $L$ act  on irreducible
graded modules $V$ as scalars, we shall use the same symbols to
denote these scalars.

\vspace{2mm} Now we  study the structure and representations of
$L_0$. Note that  by the theory of Verma modules, the irreducible
${\bf Z}$-graded highest (or lowest) weight $L$-modules are
classified by the characters of $L_0$. \vspace{2mm}

{\em {\bf Lemma 2.3 } $(1)$ If $m_{21}$ is an even integer then
$L_0$ is a Heisenberg Lie algebra.

$(2)$ If $m_{21}$ is an odd integer then
$$
L_0=({\cal A}+{\cal B})\oplus\la m_{11}c_1+m_{12}c_2\ra ,
$$
 where ${\cal
A}=\la t_0^{\bar 0}t^{2j{\bf m}_2},m_{21}c_1+m_{22}c_2\mid j\in{\bf
Z}\ra $ is a Heisenberg Lie algebra and
$${\cal B}=\la t_0^{\bar 1}t^{j{\bf
m}_2},t_0^{\bar 0}t^{(2j+1){\bf m}_2},m_{21}c_1+m_{22}c_2\mid
j\in{\bf Z}\ra ,
$$
which is isomorphic to the affine Lie algebra $A_1^{(1)}$ and the
isomorphism is given by the linear extension of the following map
$\phi$: \setcounter{section}{2}\setcounter{equation}{0}
\begin{eqnarray}
t_0^{\bar 1}t^{2j{\bf m}_2}&\mapsto&
-q^{-{2}j^2m_{22}m_{21}}((E_{11}-E_{22})(x^j)+\frac{1}{2}K),
\\%\eqno{(2.1)}$$$$
t_0^{i}t^{(2j+1){\bf m}_2}&\mapsto&
q^{-\frac{1}{2}(2j+1)^2m_{22}m_{21}}((-1)^iE_{12}(x^j)+E_{21}(x^{j+1})),
\\%\eqno{(2.1)}$$$$
m_{21}c_1+m_{22}c_2&\mapsto& K.
% \eqno{(2.3)}$$
\end{eqnarray}
Moreover, we have $[{\cal A},{\cal B}]=0$. }\vskip2mm

{\bf Proof } Statement (1) can be easily deduced from the definition
of $L_0$.

(2) To show  ${\cal B}\cong A_1^{(1)}$, we need to prove that $\phi$
preserves Lie bracket. Notice that
$$
\begin{array}{ll}
\Big[q^{-\frac{1}{2}(2j+1)^2m_{22}m_{21}}\Big((-1)^iE_{12}(x^j)+E_{21}(x^{j+1})\Big),
q^{-\frac{1}{2}(2l+1)^2m_{22}m_{21}}\Big((-1)^kE_{12}(x^l)+E_{21}(x^{l+1})\Big)\Big]
\\[7pt]
=q^{-\frac{1}{2}((2j+1)^2+(2l+1)^2)m_{22}m_{21}}\Big(
((-1)^i-(-1)^k)(E_{11}-E_{22})(x^{j+l+1})\\[7pt]
\phantom{=q^{-\frac{1}{2}((2j+1)^2+(2l+1)^2)m_{22}m_{21}}\Big(}
+\delta_{j+l+1,0}((-1)^ij+(-1)^k(j+1))K \Big), \end{array}$$ and
$$\begin{array}{ll}
[t_0^{i}t^{(2j+1){\bf m}_2},t_0^{k}t^{(2k+1){\bf m}_2}]
=\!\!\!\!&((-1)^k-(-1)^i)q^{(2j+1)(2k+1)m_{22}m_{21}}t_0^{i+k}t^{(2j+2k+2){\bf
m}_2}
\\[7pt]&
+\delta_{i+k,\bar
0}\delta_{j+k+1,0}(-1)^kq^{(2j+1)(2k+1)m_{22}m_{21}}(2j+1)(m_{21}c_1+m_{22}c_2).
\end{array}$$ One sees that $$\phi([t_0^{i}t^{(2j+1){\bf
m}_2},t_0^{k}t^{(2k+1){\bf m}_2}])=[\phi(t_0^{i}t^{(2j+1){\bf
m}_2}),\phi(t_0^{k}t^{(2k+1){\bf m}_2})].$$ Consider
$$\begin{array}{ll}
[-q^{-{2}j^2m_{22}m_{21}}((E_{11}-E_{22})(x^j)+\frac{1}{2}K),
q^{-\frac{1}{2}(2l+1)^2m_{22}m_{21}}((-1)^kE_{12}(x^l)+E_{21}(x^{l+1}))]
\\[7pt]
=-q^{-\frac{1}{2}(4j^2+(2l+1)^2)m_{22}m_{21}}(2(-1)^kE_{12}(x^{l+j})-2E_{21}(x^{l+j+1}))
\end{array}$$ and
$$
[t_0^{\bar 1}t^{2j{\bf m}_2},t_0^{k}t^{(2l+1){\bf m}_2}]=
2q^{2j(2l+1)m_{22}m_{21}}t_0^{k+\bar 1}t^{(2j+2l+1){\bf m}_2},
$$
we have $$\phi([t_0^{\bar 1}t^{2j{\bf m}_2},t_0^{k}t^{(2l+1){\bf
m}_2}])=[\phi(t_0^{\bar 1}t^{2j{\bf m}_2}),\phi(t_0^{k}t^{(2l+1){\bf
m}_2})].$$ Finally, we have
$$\begin{array}{ll}
[-q^{-{2}j^2m_{22}m_{21}}((E_{11}-E_{22})(x^j)+\frac{1}{2}K),-q^{-{2}l^2m_{22}m_{21}}((E_{11}-E_{22})(x^l)+\frac{1}{2}K)]
\\[7pt]
=2jq^{-{2}(j^2+l^2)m_{22}m_{21}}\delta_{j+l,0}K=2jq^{4jlm_{22}m_{21}}\delta_{j+l,0}K,
\end{array}$$ and
$$
[t_0^{\bar 1}t^{2j{\bf m}_2},t_0^{\bar 1}t^{2l{\bf
m}_2}]=2jq^{4jlm_{22}m_{21}}\delta_{j+l,0}(m_{21}c_1+m_{22}c_2).
$$
Thus $$\phi([t_0^{\bar 1}t^{2j{\bf m}_2},t_0^{\bar 1}t^{2l{\bf
m}_2}])=[\phi(t_0^{\bar 1}t^{2j{\bf m}_2}),\phi(t_0^{\bar
1}t^{2l{\bf m}_2})].$$ This proves ${\cal B}\cong A_1^{(1)}$. And
the proof of the rest results in this lemma is straightforward.
\hfill$\Box$

\vspace{2mm} Since the Lie subalgebra ${\cal B}$ of $L_0$ is
isomorphic to the affine Lie algebra $A_1^{(1)}$, we need to collect
some results on the finite dimensional irreducible modules of
$A_1^{(1)}$ from [30].

\def\ulai{\mu}
\def\ula{\underline{\mu}}Let $\nu>0$ and $\ula
=(\ulai_1,\ulai_2,\cdots,\ulai_{\nu})$ be a finite sequence of
nonzero distinct numbers. Let $V_i,\ 1\leq i\leq \nu$ be finite
dimensional irreducible $sl_2$-modules. We define an
$A_1^{(1)}$-module $V(\ula )=V_1\otimes V_2\otimes\cdots\otimes
V_{\nu}$ as follows, for $X\in sl_2, j\in{\bf Z}$,
$$ X(x^j)\cdot(v_1\otimes v_2\otimes\cdots\otimes
v_{\nu})=\sum\limits_{i=1}^{\nu}\ulai_i^jv_1\otimes\cdots
\otimes(X\cdot v_i)\otimes \cdots  \otimes v_{\nu},\quad
K\cdot(v_1\otimes\cdots\otimes v_{\nu})=0.
$$
Clearly $V(\ula )$ is a finite dimensional irreducible
$A_1^{(1)}$-module. For any $Q(x)\in{\bf C}[x^{\pm 1}]$, we have
$$X(Q(x))\cdot(V_1\otimes\cdots\otimes V_{\nu})=0,\;\forall\ X\in
sl_2\ \ \Longleftrightarrow\ \ \prod_{i=1}^{\nu}(x-\ulai_1)\mid
Q(x).$$ Now by Lemma 2.3(2), if $m_{21}$ is an odd integer then we
can define a finite dimensional irreducible $L_0$-module $V(\ula
,\psi)=V_1\otimes\cdots\otimes V_{\nu}$ as follows
$$\begin{array}{ll}
t_0^{\bar 0}t^{2j{\bf m}_2}\cdot (v_1\otimes\cdots\otimes
v_{\nu})=\psi(t_0^{\bar 0}t^{2j{\bf m}_2})\cdot
(v_1\otimes\cdots\otimes v_{\nu}),
\\[7pt]
t_0^{\bar 1}t^{2j{\bf m}_2}\cdot(v_1\otimes\cdots\otimes v_{\nu})=
-q^{-{2}j^2m_{22}m_{21}}\sum\limits_{i=1}^{\nu}\ulai_i^jv_1\otimes\cdots\otimes
((E_{11}-E_{22})\cdot v_i)\otimes \cdots\otimes v_{\nu},
\\[7pt]
t_0^{i}t^{(2j+1){\bf m}_2}\cdot(v_1\otimes\cdots\otimes v_{\nu})=
q^{-\frac{1}{2}(2j+1)^2m_{22}m_{21}}\Big(
(-1)^i\sum\limits_{i=1}^{\nu}\ulai_i^jv_1\otimes\cdots\otimes
(E_{12}\cdot v_i)\otimes \cdots\otimes v_{\nu}
\\[7pt]\phantom{t_0^{i}t^{(2j+1){\bf m}_2}\cdot(v_1\otimes\cdots\otimes
v_{\nu})=} +\sum\limits_{i=1}^{\nu}\ulai_i^{j+1}v_1\otimes\cdots
\otimes(E_{21}\cdot v_i)\otimes \cdots\otimes v_{\nu}\Big ),
\\[7pt]
(m_{21}c_1+m_{22}c_2)\cdot(v_1\otimes\cdots\otimes v_{\nu})=0,\quad
\forall\ v_1\otimes\cdots\otimes v_{\nu}\in V(\ula ,\psi),j\in{\bf
Z},i\in{\bf Z}_2, \end{array}$$ where $\nu>0$, $\ula
=(\ulai_1,\ulai_2,\cdots,\ulai_{\nu})$ is a finite sequence of
nonzero distinct numbers, $V_i,\ 1\leq i\leq \nu$ are finite
dimensional irreducible $sl_2$-modules, and $\psi$ is a linear
function over ${\cal A}$. \vspace{2mm}

{\em {\bf Theorem 2.4~~([30, Theorem 2.14])}~~Let $V$ be a finite
dimensional irreducible $A_1^{(1)}$-module. Then $V$ is isomorphic
to $V(\ula )$ for some finite dimensional irreducible $sl_2$-modules
$V_1,\cdots,V_{\nu}$ and a finite sequence of nonzero distinct
numbers $\ula =(\ulai_1,\cdots,\ulai_{\nu})$.} \vspace{2mm}

 From the above theorem and Lemma 2.3, we have the
following theorem.\vspace{2mm}

{\em {\bf Theorem 2.5 } Let $m_{21}$ be an odd integer and $V$ be a
finite dimensional irreducible $L_0$-module. Then $V$ is isomorphic
to $V(\ula ,\psi)$, where $V_1,\cdots,V_{\nu}$ are some finite
dimensional irreducible $sl_2$-modules, $\ula
=(\ulai_1,\cdots,\ulai_{\nu})$ is a finite sequence of nonzero
distinct numbers, and $\psi$ is a linear function over ${\cal A}$.}
\vspace{2mm}

{ {\bf Remark 2.6 }  Let $m_{21}$ be an odd integer and $V(\ula
,\psi)$ be a finite dimensional irreducible $L_0$-modules defined as
above. One can  see that for any $k\in{\bf Z}_2$,
$$\begin{array}{ll}
(\sum\limits_{i=1}^nb_iq^{\frac{1}{2}(2i+1)^2m_{22}m_{21}}t_0^{k}t^{(2i+1){\bf
m}_2})\cdot(V_1\otimes\cdots\otimes V_{\nu})=0, \mbox{ \
 and}\\[11pt]
(\sum\limits_{i=1}^nb_iq^{2i^2m_{22}m_{21}}t_0^{\bar 1}t^{2i{\bf
m}_2})\cdot(V_1\otimes\cdots\otimes V_{\nu})=0, \end{array}$$ if and
only if $\prod_{i=1}^{\nu}(x-\ulai_1)\mid (\sum_{i=1}^nb_ix^i)$.}
\vspace{2mm}

At the end of this section, we will prove a proposition which
reduces the classification of the irreducible {\bf Z}-graded modules
with finite dimensional homogeneous subspaces to that of the
generalized highest weight modules and the uniformly bounded
modules. \vspace{2mm}

{\em {\bf Proposition 2.7 } If $V$ is an irreducible {\bf Z}-graded
$L$-module, then $V$ is a generalized highest weight module or a
uniformly bounded module.} \vspace{2mm}

{\bf Proof } Let $V=\oplus_{m\in{\bf Z}}V_m$. We first prove that if
there exists a  ${\bf Z}$-basis  $\{{\bf b}_1, {\bf b}_2\}$ of ${\bf
Z}^2$ and a homogeneous vector $v\neq 0$
 such that
 $t_0^it^{{\bf b}_1}\cdot v=t_0^it^{{\bf b}_2}\cdot v=0,\;\forall\ i\in{\bf Z}/2{\bf Z}$, then $V$ is a
generalized highest weight modules.

To obtain this, we first introduce the following notation: For any
$A\subset{\bf Z}^2$, we use $t^A$ to denote the set $\{t^{a}|{ a}\in
A\}$.

Then one can deduce that $t_0^it^{{\bf Nb}_1+{\bf Nb}_2}\cdot
v=0,\;\forall\ i\in{\bf Z}/2{\bf Z}$ by the assumption. Thus for the
${\bf Z}$-basis ${\bf m}_1=3{\bf b}_1+{\bf b}_2,\; {\bf m}_2=2{\bf
b}_1+{\bf b}_2$ of ${\bf Z}^2$ we have $t_0^it^{{\bf Z}_+{\bf
m}_1+{\bf Z}_+{\bf m}_2}v=0,\;\forall\ i\in{\bf Z}_2$. Therefore $V$
is a generalized highest weight module by the definition.

With the above statement, we can prove our proposition now. Suppose
that $V$ is not a generalized highest weight module. For any
$m\in{\bf Z},$ considering the following maps
$$\begin{array}{llllll}
t_0^{{\bar 0}}t^{-m{\bf m}_1+{\bf m}_2}:&V_m\mapsto V_0,\ \ \ \ \ \
\ \ \ \ \quad& t_0^{{\bar 1}}t^{-m{\bf m}_1+{\bf m}_2}:&V_m\mapsto
V_0,
\\[7pt]
t_0^{{\bar 0}}t^{(1-m){\bf m}_1+{\bf m}_2}:&V_m\mapsto
V_1,\quad&t_0^{{\bar 1}}t^{(1-m){\bf m}_1+{\bf m}_2}:&V_m\mapsto
V_1, \end{array}$$ one can easily check that
$$
\mbox{ker\,}t_0^{{\bar 0}}t^{-m{\bf m}_1+{\bf
m}_2}\cap\mbox{ker\,}t_0^{{\bar 0}}t^{(1-m){\bf m}_1+{\bf m}_2}\cap
\mbox{ker\,}t_0^{{\bar 1}}t^{-m{\bf m}_1+{\bf
m}_2}\cap\mbox{ker\,}t_0^{{\bar 1}}t^{(1-m){\bf m}_1+{\bf
m}_2}=\{0\}.
$$
Therefore $\mbox{dim}V_m\leq 2\mbox{dim}V_0+2\mbox{dim}V_1$. So $V$
is a uniformly bounded module.\hfill$\Box$
 \vspace{7pt}

\centerline{\large\bf \S 3 \ \ The highest weight irreducible ${\bf
Z}$-graded $L$-modules }

For any finite dimensional irreducible ${ L}_{0}$-module $V$,
 we can define it as a
$({L}_{0}+{L}_{+})$-module by putting $L_+v=0,\;\forall\ v\in V$.
Then we obtain an induced ${L}$-module, \def\ol{\overline}$$ \ol
M{}^{+}(V,{\bf m}_1,{\bf m}_2)
=\mbox{Ind}^{{L}}_{{L}_{0}+{L}_{+}}V=U({{L}})\otimes_{U({L}_{0}+{L}_{+})}V\simeq
U({L}_{-})\otimes V,
$$
where $U({L})$ is the universal enveloping algebra of ${L}$. If we
set $V$ to be the homogeneous subspace of $\ol M{}^{+}(V,{\bf
m}_1,{\bf m}_2)$ with degree $0$, then $\ol M{}^{+}(V,{\bf m}_1,{\bf
m}_2)$ becomes a ${\bf Z}$-graded $L$-module in a natural way.
Obviously,
 $\ol
M{}^{+}(V,{\bf m}_1,{\bf m}_2)$ has an unique maximal proper
submodule $J$ which trivially intersects with $V$. So we obtain an
irreducible {\bf Z}-graded highest weight
${L}$-module,\def\OVERLINE{}
$$
\OVERLINE{M}^{+}(V,{\bf m}_1,{\bf m}_2)=\ol M{}^{+}(V,{\bf m}_1,{\bf
m}_2)/J.
$$
We can write it as
$$
\OVERLINE{M}^{+}(V,{\bf m}_1,{\bf m}_2)=\bigoplus_{i\in{\bf
Z}_+}V_{-i},
$$
where $V_{-i}$ is the homogeneous subspaces of degree $-i$. Since
$L_-$ is generated by $L_{-1}$, and $L_{+}$ is generated by $L_{1}$,
by the construction of $ \OVERLINE{M}^{+}(V,{\bf m}_1,{\bf m}_2)$,
we see that
$$
L_{-1}V_{-i}=V_{-i-1},\quad \forall \ i\in{\bf Z}_+,
\eqno{(3.1)}
$$
and for a homogeneous vector $v$,
$$
 L_1\cdot v=0\ \Longrightarrow\ v=0.
\eqno{(3.2)}
$$

Similarly, we can define an irreducible lowest weight ${\bf
Z}$-graded ${L}$-module $\OVERLINE{M}^{-}(V,{\bf m}_1,{\bf m}_2)$
from any finite dimensional irreducible ${L}_{0}$-module $V$.

If $m_{21}\in{\bf Z}$ is even then $L_0$ is a Heisenberg Lie algebra
by Lemma 2.3. Thus, by a well-known result about the representations
of the Heisenberg Lie algebra, we see that the finite dimensional
irreducible $L_0$-module $V$ must be a one dimensional module ${\bf
C}v_0$, and there is a linear function $\psi$ over $L_0$ such that
$$
t_0^it^{j{\bf m}_2}\cdot v_0=\psi(t_0^it^{j{\bf m}_2})\cdot v_0,\;\;
\psi(m_{21}c_1+m_{22}c_2)=0,\,\,\forall\ i\in{\bf Z}_2, j\in{\bf Z}.
$$
In this case, we denote the corresponding highest weight, resp.,
lowest weight, irreducible {\bf Z}-graded $L$-module by
$$\OVERLINE{M}^+(\psi,{\bf m}_1,{\bf m}_2), \mbox{ \ \ \ resp., \ \ \ }
\OVERLINE{M}^-(\psi,{\bf m}_1,{\bf m}_2).$$

If $m_{21}$ is an odd integer then $V$ must be isomorphic to $V(\ula
,\psi)$ by Theorem 2.5. We denote the corresponding highest weight,
resp.~lowest weight, irreducible {\bf Z}-graded $L$-module by
$$\OVERLINE{M}^+(\ula ,{\psi},{\bf m}_1,{\bf m}_2),\mbox{ \
\ \ resp., \ \ \ }\OVERLINE{M}^-(\ula ,{\psi},{\bf m}_1,{\bf
m}_2).$$

The irreducible {\bf Z}-graded $L$-modules $\OVERLINE{M}^+(\psi,{\bf
m}_1,{\bf m}_2)$ and $\OVERLINE{M}^+(\ula ,{\psi},{\bf m}_1,{\bf
m}_2)$ are in general not quasi-finite modules. Thus in the rest of
this section we shall determine which of $\ula $ and $\psi$ can
correspond to quasi-finite modules.

For the later use, we obtain the following equations from the
definition of $L$, where, $\alpha=m_{11}m_{22}-m_{12}m_{21}\in\{\pm
1\}$, \setcounter{section}{3}\setcounter{equation}{2}
\begin{eqnarray}&&
[t_0^jt^{{\bf m}_1+k{\bf m}_2},t_0^rt^{-{\bf m}_1+s{\bf
m}_2}t^{i{\bf m}_2}] \nonumber\\&&
=q^{i(-m_{12}+sm_{22})m_{21}}[t_0^jt^{{\bf
m}_1+k{\bf m}_2},t_0^rt^{-{\bf m}_1+(s+i){\bf m}_2}]\nonumber \\
&&=q^{-m_{11}m_{12}-km_{11}m_{22}+sm_{12}m_{21}+ksm_{21}m_{22}}(-1)^{r(m_{11}+km_{21})}\times
\nonumber\\
 &&\phantom{=} \times\left((1-(-1)^{(j+r)m_{11}+(kr+js+ji)m_{21}}q^{(k+s+i)\alpha})t_0^{j+r}t^{(k+s){\bf
m}_2}t^{i{\bf m}_2}
\right.\nonumber\\
&&\phantom{=} \left.+\delta_{k+s+i,0}\delta_{j+r,\bar
0}q^{-(k+s)^2m_{21}m_{22}}((m_{11}+km_{21})c_1+(m_{12}+km_{22})c_2)\right),
\end{eqnarray}\begin{eqnarray}%\\[10pt]
&&
 [t_0^st^{k{\bf m}_2}t^{i{\bf m}_2},t_0^rt^{-{\bf m}_1+j{\bf
m}_2}]\nonumber\\
&&=q^{kim_{22}m_{21}}[t_0^st^{(k+i){\bf m}_2},t_0^rt^{-{\bf
m}_1+j{\bf m}_2}]
\nonumber\\
&&=q^{km_{22}(-m_{11}+jm_{21})}(-1)^{(rk+ri)m_{21}}(q^{-i\alpha}-(-1)^{sm_{11}+(rk+ri+sj)m_{21}}q^{k\alpha})
\times \nonumber\\&&\phantom{=}\times t_0^{r+s}t^{-{\bf
m}_1+(k+j){\bf m}_2}t^{i{\bf m}_2}.
%\eqno{(3.4)}
\end{eqnarray}

In the rest of this section, if $P(t^{{\bf
m}_2})=\sum_{i=0}^{n}a_it^{i{\bf m}_2}\in{\bf C}[t^{{\bf m}_2}]$
then we will denote $\sum_{i=0}^{n}a_ib^it^{i{\bf m}_2}$ by
$P(bt^{{\bf m}_2})$ for any $b\in{\bf C}$.\vspace*{2mm}

{\em {\bf Lemma 3.1 } Let $m_{21}$ be an even integer. Then
$\OVERLINE{M}^\pm(\psi,{\bf m}_1,{\bf m}_2)
%,\OVERLINE{M}^-(\psi,{\bf m}_1,{\bf m}_2)$ $
\in{\cal O}_{\bf Z}$ if and only if there exists a polynomial
$P(t^{{\bf m}_2})=\sum_{i=0}^{n}a_it^{i{\bf m}_2}\in{\bf C}[t^{{\bf
m}_2}]$ with $a_0a_n\neq 0$ such that for  $k\in{\bf Z}, j\in{\bf
Z}_2$,
$$
\psi\Big(t_0^jt^{k{\bf m}_2}P(t^{{\bf
m}_2})-(-1)^{j}q^{k\alpha}t_0^jt^{k{\bf m}_2}P(q^{\alpha}t^{{\bf
m}_2})+\delta_{j,{\bar 0}}a_{-k}q^{-k^2m_{21}m_{22}}\beta\Big)=0,
\eqno{(3.5)}
$$
where $a_k=0$ if $k\not\in\{0,1,\cdots,n\}$, and
$\alpha=m_{11}m_{22}-m_{12}m_{21}\in\{\pm 1\}$, $\beta=
m_{11}c_1+m_{12}c_2$.}\vspace*{2mm}

{\bf Proof}\quad Since $m_{21}$ is an even integer and
$m_{11}m_{22}-m_{12}m_{21}\in\{\pm 1\}$, we see $m_{11}$ is an odd
integer.

``$\Longrightarrow$''.\quad Since $\mbox{dim}V_{-1}<\infty$, there
exist an integer $s$ and a polynomial $P(t^{{\bf
m}_2})=\sum_{i=0}^{n}a_it^{i{\bf m}_2}$ $\in{\bf C}[t^{{\bf m}_2}]$
with $a_0a_n\neq0$ such that
$$
t_0^{\bar{0}}t^{-{\bf m}_1+s{\bf m}_2}P(t^{{\bf m}_2})\cdot v_0=0.
$$
Applying $t_0^jt^{{\bf m}_1+k{\bf m}_2}$ for any $k\in{\bf Z},
j\in{\bf Z}_2$ to the above equation, we have that
$$
0=t_0^jt^{{\bf m}_1+k{\bf m}_2}\cdot t_0^{{\bar 0}}t^{-{\bf
m}_1+s{\bf m}_2}P(t^{{\bf m}_2})\cdot v_0
=\sum\limits_{i=0}^{n}[t_0^jt^{{\bf m}_1+k{\bf m}_2},a_it_0^{{\bar
0}}t^{-{\bf m}_1+s{\bf m}_2}t^{i{\bf m}_2}]\cdot v_0.
$$
Thus, by (3.3), we have
$$\begin{array}{ll}
0\!\!\!\!&=
\psi\left(\sum\limits_{i=0}^{n}a_i\left((1-(-1)^{j}q^{(k+s+i)\alpha})t_0^{j}t^{(k+s){\bf
m}_2}t^{i{\bf m}_2}+\delta_{k+s+i,0}\delta_{j,\bar
0}q^{-(k+s)^2m_{21}m_{22}}\beta\right)\right)
\\[12pt]&
=\psi\Big(t_0^{j}t^{(k+s){\bf m}_2}P(t^{{\bf
m}_2})-(-1)^{j}q^{(k+s)\alpha}t_0^{j}t^{(k+s){\bf
m}_2}P(q^{\alpha}t^{{\bf m}_2})+a_{-k-s}\delta_{j,\bar
0}q^{-(k+s)^2m_{21}m_{22}}\beta\Big). \end{array}$$ Therefore this
direction follows.\vskip2mm

``$\Longleftarrow$''.\quad By induction on $s$ we first show the
following claim.\vskip2mm

{\em {\bf Claim.}\quad For any $s\in{\bf Z}_+$, there exists
polynomial $P_s(t^{{\bf m}_2})=\sum_{i\in{\bf Z}}a_{s,i}t^{i{\bf
m}_2}\in{\bf C}[t^{{\bf m}_2}]$ such that
$$\begin{array}{ll}
\Big(t_0^rt^{k{\bf m}_2}P_s(t^{{\bf
m}_2})-(-1)^{r}q^{k\alpha}t_0^rt^{k{\bf m}_2}P_s(q^{\alpha}t^{{\bf
m}_2})+\delta_{r,{\bar
0}}a_{s,-k}q^{-k^2m_{21}m_{22}}\beta\Big)\cdot V_{-s}=0,
\\[7pt]
t_0^rt^{-{\bf m}_1+k{\bf m}_2}P_s(t^{{\bf m}_2})\cdot V_{-s}=0,\;\;\
\ \forall\ r\in{\bf Z}_2,k\in{\bf Z}. \end{array}$$ }

For $s=0$, the first equation holds with $P_0(t^{{\bf
m}_2})=P(t^{{\bf m}_2})$ (with $P$ being as in the necessity), and
by (3.2),  the second equation can be deduced by a calculation
similar to the proof of the necessity. Suppose the claim holds for
$s$. Let us consider the claim for $s+1$.

Note that the equations in the claim are equivalent to
\setcounter{equation}{5}\begin{eqnarray}&& \Big(t_0^rQ(t^{{\bf
m}_2})-(-1)^{r}t_0^rQ(q^{\alpha}t^{{\bf m}_2})+\delta_{r,{\bar
0}}a_{Q}\beta\Big)\cdot V_{-s}=0, \nonumber\\&& t_0^rt^{-{\bf
m}_1+k{\bf m}_2}Q(t^{{\bf m}_2})\cdot V_{-s}=0,\;\;\ \forall\
r\in{\bf Z}_2,k\in{\bf Z},
\end{eqnarray}%\eqno{(3.6)}$$
for any $Q(t^{{\bf m}_2})\in{\bf C}[t^{\pm{\bf m}_2}]$ with
 $P_s(t^{{\bf m}_2})\mid Q(t^{{\bf
m}_2})$, where $a_{Q}$ is the constant term of $Q(t^{{\bf m}_2})$.

Let $P_{s+1}(t^{{\bf m}_2})=P_s(q^{\alpha}t^{{\bf m}_2})P_s(t^{{\bf
m}_2})P_s(q^{-\alpha}t^{{\bf m}_2})$, then $$\mbox{$P_s(t^{{\bf
m}_2})\mid P_{s+1}(t^{{\bf m}_2}),\ \ \   P_s(t^{{\bf m}_2})\mid
P_{s+1}(q^{\alpha}t^{{\bf m}_2})$ \ \ and \ \ $P_s(t^{{\bf
m}_2})\mid P_{s+1}(q^{-\alpha}t^{{\bf m}_2})$.}$$ For any
$p,r\in{\bf Z}_2,\, j,k\in{\bf Z}$, by induction and (3.4), we have
$$\begin{array}{lllll}
\Big(t_0^rt^{k{\bf m}_2}P_{s+1}(t^{{\bf
m}_2})-(-1)^{r}q^{k\alpha}t_0^rt^{k{\bf
m}_2}P_{s+1}(q^{\alpha}t^{{\bf m}_2})+\delta_{r,{\bar
0}}a_{s+1,-k}q^{-k^2m_{21}m_{22}}\beta\Big)\cdot t_0^pt^{-{\bf
m}_1+j{\bf m}_2}\cdot V_{-s}
\\[7pt]
=\Big[t_0^rt^{k{\bf m}_2}P_{s+1}(t^{{\bf
m}_2})-(-1)^{r}q^{k\alpha}t_0^rt^{k{\bf
m}_2}P_{s+1}(q^{\alpha}t^{{\bf m}_2})+\delta_{r,{\bar
0}}a_{s+1,-k}q^{-k^2m_{21}m_{22}}\beta,t_0^pt^{-{\bf m}_1+j{\bf
m}_2}\Big]\cdot V_{-s}
\\[7pt]=q^{-km_{22}m_{11}+kjm_{22}m_{21}}\left(t_0^{r+p}t^{-{\bf
m}_1+(k+j){\bf m}_2}\Big(P_{s+1}(q^{-\alpha}t^{{\bf
m}_2})-2(-1)^rq^{k\alpha}P_{s+1}(t^{{\bf m}_2}) \right.
\\[7pt]
%\phantom{=q^{-km_{22}m_{11}+kjm_{22}m_{21}}\Big(t_0^{r+p}t^{-{\bf
%m}_1+(k+j){\bf m}_2}\left(\right.}
\phantom{=}\left. +q^{2k\alpha}P_{s+1}(q^{\alpha}t^{{\bf
m}_2})\Big)\right)\cdot
V_{-s}\\[7pt]=0. \end{array}$$ Thus, by (3.1) and (3.2), we obtain that
$$
\Big(t_0^rt^{k{\bf m}_2}P_{s+1}(t^{{\bf
m}_2})-(-1)^{r}q^{k\alpha}t_0^rt^{k{\bf
m}_2}P_{s+1}(q^{\alpha}t^{{\bf m}_2})+\delta_{r,{\bar
0}}a_{s+1,-k}q^{-k^2m_{21}m_{22}}\beta\Big)\cdot V_{-s-1}=0.
\eqno{(3.7)}
$$
This proves the first equation in the claim for $i=s+1$.

Using (3.3), (3.6) and induction, we deduce that for any $l,k\in{\bf
Z},\,\,n,r\in{\bf Z}_2$,
\begin{eqnarray*}&&
t_0^nt^{{\bf m}_1+l{\bf m}_2}\cdot t_0^rt^{-{\bf m}_1+k{\bf
m}_2}P_{s+1}(t^{{\bf m}_2})\cdot   V_{-s-1}
\\[7pt]&&
=[t_0^nt^{{\bf m}_1+l{\bf m}_2},t_0^rt^{-{\bf m}_1+k{\bf
m}_2}P_{s+1}(t^{{\bf m}_2})]\cdot V_{-s-1}
%\\&& \phantom{=}
+t_0^rt^{-{\bf m}_1+k{\bf
m}_2}P_{s+1}(t^{{\bf m}_2})\cdot t_0^nt^{{\bf m}_1+l{\bf m}_2}\cdot
 V_{-s-1}
\\[7pt]&&
=(-1)^rq^{-m_{11}m_{12}+km_{12}m_{21}-lm_{11}m_{22}+lkm_{21}m_{22}}\Big(t_0^{n+r}t^{(l+k){\bf
m}_2}P_{s+1}(t^{{\bf m}_2})
\\&&\phantom{=}-
(-1)^{n+r}q^{(k+l)\alpha}t_0^{n+r}t^{(l+k){\bf
m}_2}P_{s+1}(q^{\alpha}t^{{\bf m}_2})+a_{s+1,-l-k}\delta_{r+n,{\bar
0}}q^{-(l+k)^2m_{21}m_{22}}\beta\Big)\cdot V_{-s-1}\\[7pt]&&=0,
\end{eqnarray*}
since $t_0^nt^{{\bf m}_1+l{\bf m}_2}\cdot V_{-s-1}\in V_{-s} $.
Hence by (3.2),$$\mbox{$t_0^rt^{-{\bf m}_1+k{\bf
m}_2}P_{s+1}(t^{{\bf m}_2})\cdot V_{-s-1}=0$ for all $r\in{\bf
Z}_2,\;k\in{\bf Z}$,}$$  which implies the second equation in the
claim for $i=s+1$. Therefore the claim follows by induction.

From the second equation of the claim and (3.1), we see that
$$
\mbox{dim}V_{-s-1}\leq 2\mbox{deg}(P_{s+1}(t^{{\bf m}_2}))\cdot
\mbox{dim}V_{s},\;\;\forall \ s\in{\bf Z}_+,
$$
where $\mbox{deg}(P_{s+1}(t^{{\bf m}_2}))$ denotes the degree of
polynomial $P_{s+1}(t^{{\bf m}_2})$. Hence $\OVERLINE{M}^+(\psi,{\bf
m}_1,{\bf m}_2)\in{\cal O}_{\bf Z}$.

Similarly we can prove the statement for $\OVERLINE{M}^-(\psi,{\bf
m}_1,{\bf m}_2)$. \hfill$\Box$\vskip2mm

{\em{\bf Theorem 3.2}\quad Let $m_{21}$ be an even integer. Then
$\OVERLINE{M}^\pm(\psi,{\bf m}_1,{\bf m}_2)\in{\cal O}_{\bf Z}$ if
and only if there exist
$b_{10}^{(j)},b_{11}^{(j)},\cdots,b_{1s_1}^{(j)},b_{20}^{(j)},b_{21}^{(j)},\cdots,b_{2s_2}^{(j)},\cdots,b_{r0}^{(j)},b_{r1}^{(j)},\cdots,b_{rs_r}^{(j)}\in{\bf
C}\,$ for $j\in{\bf Z}_2$, and $\alpha_1,\cdots,\alpha_r\in{\bf
C}^{*}$ such that for any $i\in{\bf Z}^{*}$, $j\in{\bf Z}_2$,
$$\begin{array}{lll}\displaystyle
\psi(t_0^jt^{i{\bf
m}_2})=\frac{(b_{10}^{(j)}+b_{11}^{(j)}i+\cdots+b_{1s_1}^{(j)}i^{s_1})\alpha_1^i
+\cdots+(b_{r0}^{(j)}+b_{r1}^{(j)}i+\cdots,b_{rs_r}^{(j)}i^{s_r})\alpha_r^i}{(1-(-1)^jq^{i\alpha
})q^{\frac{1}{2}i^2m_{21}m_{22}}},
\\[7pt]\psi(\beta)=b_{10}^{(0)}+b_{20}^{(0)}+\cdots+b_{r0}^{(0)},
\\[7pt]\psi(t_0^{\bar 1}t^{\bf
0})=\frac{1}{2}(b_{10}^{(1)}++b_{20}^{(1)}+\cdots+b_{r0}^{(1)}),
\mbox{ \ \ \ and \ \ \ }\psi(m_{21}c_1+m_{22}c_2)=0,
\end{array}$$
 where $\alpha=m_{11}m_{22}-m_{21}m_{12}\in\{\pm
1\}$ and $\beta=m_{11}c_1+m_{12}c_2$.}\vskip2mm

{\bf Proof}\quad ``$\Longrightarrow$''.\quad Let
$f_{j,i}=\psi((1-(-1)^jq^{i\alpha
})q^{\frac{1}{2}i^2m_{21}m_{22}}t_0^jt^{i{\bf m}_2})$ for $j\in{\bf
Z}_2,\;i\in{\bf Z}^{*}$   and $f_{0,0}=\psi(\beta),\;f_{1,0}=\psi
(2t_0^{1^{}}t^{\bf 0^{^{}}})$. By Lemma 3.1 there exist complex
numbers $a_0,a_1,\cdots,a_n$ with $a_0a_n\neq 0$ such that
$$
\sum_{i=0}^na_iq^{-\frac{1}{2}i^2m_{21}m_{22}}f_{j,k+i}=0,\;\;\forall\
k\in{\bf Z},j\in{\bf Z}_2.
 \eqno{(3.8)}
$$
Denote $b_i=a_iq^{-\frac{1}{2}i^2m_{21}m_{22}}$. Then the above
equation becomes
$$
\sum_{i=0}^nb_if_{j,k+i}=0,\;\;\forall\ k\in{\bf Z},j\in{\bf Z}_2.
\eqno{(3.9)}
$$
Suppose $\alpha_1,\cdots,\alpha_r$ are all distinct roots of the
equation $\sum_{i=0}^nb_ix^i=0$ with multiplicity
$s_1+1,\cdots,s_r+1$ respectively. By a well-known combinatorial
formula, we see that there exist
$b_{10}^{(j)},b_{11}^{(j)},\cdots,b_{1s_1}^{(j)},\cdots,b_{r0}^{(j)},b_{r1}^{(j)},\cdots,b_{rs_r}^{(j)}\in{\bf
C}$ for $j\in{\bf Z}_2$ such that
$$
f_{j,i}=(b_{10}^{(j)}+b_{11}^{(j)}i+\cdots+b_{1s_1}^{(j)}i^{s_1})\alpha_1^i
+\cdots+(b_{r0}^{(j)}+b_{r1}^{(j)}i+\cdots,b_{rs_r}^{(j)}i^{s_r})\alpha_r^i,
\;\;\forall\ i\in{\bf Z}.
$$
Therefore, for any $i\in{\bf Z}^{*},\ j\in{\bf Z}_2$,
$$\begin{array}{llll}\displaystyle
\psi(t_0^jt^{i{\bf
m}_2})=\frac{(b_{10}^{(j)}+b_{11}^{(j)}i+\cdots+b_{1s_1}^{(j)}i^{s_1})\alpha_1^i
+\cdots+(b_{r0}^{(j)}+b_{r1}^{(j)}i+\cdots,b_{rs_r}^{(j)}i^{s_r})\alpha_r^i}{(1-(-1)^jq^{i\alpha
})q^{\frac{1}{2}i^2m_{21}m_{22}}},
\\[11pt]\psi(\beta)=f_{0,0}=b_{10}^{(0)}+b_{20}^{(0)}+\cdots+b_{r0}^{(0)},
\mbox{ \ \  and}
\\[7pt]
\psi(t_0^{\bar 1}t^{\bf
0})=f_{1,0}=\frac{1}{2}(b_{10}^{(1)}++b_{20}^{(1)}+\cdots+b_{r0}^{(1)}).
\end{array}$$ Thus we obtain the expression as required. This direction
follows.\vskip2mm

``$\Longleftarrow$''. Set
$$
Q(x)=\prod_{i=1}^{r}(x-\alpha_i)^{s_i+1}=\sum_{i=1}^nb_ix^i\in{\bf
C}[x], \;\; f_{j,i}=(1-(-1)^jq^{i\alpha
})q^{\frac{1}{2}i^2m_{21}m_{22}}\psi(t_0^jt^{i{\bf m}_2}),
$$
for $j\in{\bf Z}_2,\;i\in{\bf Z}^{*}$, and set
$$
f_{0,0}=\psi(\beta),\;f_{1,0}=2\psi (t_0^1t^{\bf 0}).
$$
Then one can verify that (3.9) holds. Let
$a_i=q^{\frac{1}{2}i^2m_{21}m_{22}}b_i$. One deduces that (3.8)
holds. Thus (3.5) holds for $P(t^{{\bf
m}_2})=\sum_{i=0}^na_it^{i{\bf m}_2}$. Therefore this direction
follows by using Lemma 3.1. %This completes the proof of the theorem.
\hfill$\Box$

\vspace{2mm} {\em {\bf Lemma 3.3}\quad  If $m_{21}$ is an odd
integer, then $\OVERLINE{M}^+(\ula ,\psi,{\bf m}_1,{\bf
m}_2)\in{\cal O}_{\bf Z}$ if and only if there exists a polynomial
$P(t^{{\bf m}_2})=\sum_{i=0}^{n}a_it^{2i{\bf m}_2}\in{\bf C}[t^{{\bf
m}_2}]$ with $a_0a_n\neq 0$ such that for any $k\in{\bf Z}$ and
$v\in V_0$,\setcounter{equation}{9}
\begin{eqnarray}&&
\Big(t_0^{\bar 0}t^{2k{\bf m}_2}P(t^{{\bf
m}_2})-q^{2k\alpha}t_0^{\bar 0}t^{2k{\bf m}_2}P(q^{\alpha}t^{{\bf
m}_2})+a_{-k}q^{-4k^2m_{21}m_{22}}\beta\Big)\cdot v=0,
%\eqno{(3.10)}
\\[0pt]&&
t_0^{\bar 0}t^{(2k+1){\bf m}_2}P(t^{{\bf m}_2})\cdot v=t_0^{\bar
0}t^{(2k+1){\bf m}_2}P(q^{\alpha}t^{{\bf m}_2})\cdot v=0,
\\[0pt]&&%\eqno{(3.11)}$$$$
t_0^{\bar 1}t^{k{\bf m}_2}P(t^{{\bf m}_2})\cdot v=t_0^{\bar
1}t^{k{\bf m}_2}P(q^{\alpha}t^{{\bf m}_2})\cdot v =0, %\eqno{(3.12)}
\end{eqnarray}
where $a_k=0$ if $k\not\in\{0,1,\cdots,n\}$, and
$\alpha=m_{11}m_{22}-m_{12}m_{21}$,  $\beta=
m_{11}c_1+m_{12}c_2$.}\vskip2mm

{\bf Proof}\quad ``$\Longrightarrow$''. Since $V_0$ is a finite
dimensional irreducible $L_0$-module, we have $V_0\cong V(\ula
,\psi)$ as $L_0$-modules by Theorem 2.5. Since ${\cal H}=\la
t_0^{\bar 1}t^{2k{\bf m}_2}\mid k\in{\bf Z}\ra $ is an Abelian Lie
subalgebra of $L_0$, we can choose a common eigenvector $v_0\in V_0$
of $\cal H$. First we prove the following claim.\vskip2mm

{\bf Claim 1} There is a polynomial $P_{e}(t^{{\bf
m}_2})=\sum_{i=0}^na_it^{2i{\bf m}_2}$ with $a_na_0\neq 0$ such that
\begin{eqnarray}
&&\Big(t_0^{\bar 0}t^{2k{\bf m}_2}Q(t^{{\bf m}_2})-
q^{2k\alpha}t_0^{\bar 0}t^{2k{\bf m}_2}Q(q^{\alpha}t^{{\bf m}_2})
+a_{Q}\beta\Big)\cdot v_0 =0, \nonumber\\[0pt]&&%$$$$
\Big(t_0^{\bar 1}t^{2k{\bf m}_2}Q(t^{{\bf m}_2})-(-1)^{m_{11}}
q^{2k\alpha}t_0^{\bar 1}t^{2k{\bf m}_2}Q(q^{\alpha}t^{{\bf
m}_2})\Big)\cdot v_0 =0,
\nonumber\\[0pt]&&%$$$$
\Big(t_0^{\bar 0}t^{(2k+1){\bf m}_2}Q(t^{{\bf
m}_2})-q^{(2k+1)\alpha}t_0^{\bar 0}t^{(2k+1){\bf
m}_2}Q(q^{\alpha}t^{{\bf m}_2})\Big)\cdot v_0=0, \nonumber
%\end{eqnarray}\begin{eqnarray}
\\[0pt]
&&%$$$$
\Big(t_0^{\bar 1}t^{(2k+1){\bf m}_2}Q(t^{{\bf
m}_2})-(-1)^{m_{11}}q^{(2k+1)\alpha}t_0^{\bar 0}t^{(2k+1){\bf
m}_2}Q(q^{\alpha}t^{{\bf m}_2})\Big)\cdot v_0=0,\end{eqnarray}
%\eqno{(3.13)}$$
for all $k\in{\bf Z}$ and $Q(t^{{\bf m}_2})\in{\bf
C}[t^{\pm 2{\bf m}_2}]$ with $P_{e}(t^{{\bf m}_2})\mid Q(t^{{\bf
m}_2})$, where $a_{Q}$ is the constant term of $t^{2k{\bf
m}_2}Q(t^{{\bf m}_2})$.

\vskip2mm

To prove the claim, since $\mbox{dim}V_{-1}<\infty$, there exist an
integer $s$ and a polynomial $P_{e}(t^{{\bf
m}_2})=\sum_{i=0}^{n}a_it^{2i{\bf m}_2}\in{\bf C}[t^{{\bf m}_2}]$
with $a_0a_n\neq0$ such that
$$
t_0^{\bar{0}}t^{-{\bf m}_1+2s{\bf m}_2}P_{e}(t^{{\bf m}_2})\cdot
v_0=0. \eqno{(3.14)}
$$
Applying $t_0^{\bar 0}t^{{\bf m}_1+2k{\bf m}_2}$ for any $k\in{\bf
Z}$ to the above equation, we have\setcounter{equation}{14}
\begin{eqnarray}\!\!\!\!\!\!\!\!\!\!\!\!&\!\!\!\!\!\!\!\!\!\!\!\!\!\!\!\!\!&
0=t_0^{\bar 0}t^{{\bf m}_1+2k{\bf m}_2}\cdot t_0^{{\bar 0}}t^{-{\bf
m}_1+2s{\bf m}_2}P_{e}(t^{{\bf m}_2})\cdot v_0 \nonumber\\[0pt]\!\!\!\!\!\!\!\!\!\!\!\!&\!\!\!\!\!\!\!\!\!\!\!\!\!\!\!\!\!\!\!\!\!\!\!\!&
\phantom{0}=\mbox{$\sum\limits_{i=0}^{n}$}a_i[t_0^{\bar 0}t^{{\bf
m}_1+2k{\bf m}_2},q^{2im_{21}(-m_{12}+2sm_{22})}t_0^{{\bar
0}}t^{-{\bf m}_1+2(s+i){\bf m}_2}]\cdot v_0
\nonumber\\[7pt]\!\!\!\!\!\!\!\!\!\!\!\!&\!\!\!\!\!\!\!\!\!\!\!\!\!\!\!\!\!\!\!\!\!\!\!\!&
\phantom{0}=q^{-m_{11}m_{12}-2km_{22}m_{11}+2sm_{12}m_{21}+4ksm_{21}m_{22}}\times
\nonumber\\[0pt]\!\!\!\!\!\!\!\!\!\!\!\!&\!\!\!\!\!\!\!\!\!\!\!\!\!\!\!\!\!\!\!\!\!\!\!\!&
\phantom{0=}\times\Big(t_0^{\bar 0}t^{2(k+s){\bf m}_2}P_{e}(t^{{\bf
m}_2})- q^{2(s+k)\alpha}t_0^{\bar 0}t^{2(k+s){\bf
m}_2}P_{e}(q^{\alpha}t^{{\bf m}_2})
+a_{-k-s}q^{-4(k+s)^2m_{21}m_{22}}\beta\Big)\cdot v_0.\end{eqnarray}
Now applying $t_0^{\bar 1}t^{{\bf m}_1+2k{\bf m}_2}$ for any
$k\in{\bf Z}$ to (3.14), we have
\begin{eqnarray}\!\!\!\!\!\!\!\!\!\!\!\!&\!\!\!\!\!\!\!\!\!\!\!\!\!\!\!\!\!&
0=t_0^{\bar 1}t^{{\bf m}_1+2k{\bf m}_2}\cdot t_0^{{\bar 0}}t^{-{\bf
m}_1+2s{\bf m}_2}P_{e}(t^{{\bf m}_2})\cdot v_0
\nonumber\\[7pt]\!\!\!\!\!\!\!\!\!\!\!\!&\!\!\!\!\!\!\!\!\!\!\!\!\!\!\!\!\!\!\!\!\!\!\!\!&
\phantom{0} =\mbox{$\sum\limits_{i=0}^{n}$}a_i[t_0^{\bar 1}t^{{\bf
m}_1+2k{\bf m}_2},q^{2im_{21}(-m_{12}+2sm_{22})}t_0^{{\bar
0}}t^{-{\bf m}_1+2(s+i){\bf m}_2}]\cdot v_0
\nonumber\\[7pt]\!\!\!\!\!\!\!\!\!\!\!\!&\!\!\!\!\!\!\!\!\!\!\!\!\!\!\!\!\!\!\!\!\!\!\!\!&
\phantom{0}
=q^{-m_{11}m_{12}-2km_{22}m_{11}+2sm_{12}m_{21}+4ksm_{21}m_{22}}\times
\nonumber\\[0pt]\!\!\!\!\!\!\!\!\!\!\!\!&\!\!\!\!\!\!\!\!\!\!\!\!\!\!\!\!\!\!\!\!\!\!\!\!&
\phantom{0=}\times \Big(t_0^{\bar 1}t^{2(k+s){\bf m}_2}P_{e}(t^{{\bf
m}_2})-(-1)^{m_{11}} q^{2(s+k)\alpha}t_0^{\bar 1}t^{2(k+s){\bf
m}_2}P_{e}(q^{\alpha}t^{{\bf m}_2})\Big)\cdot v_0.\end{eqnarray}
%\eqno{(3.16)}$$
By applying $t_0^{\bar 0}t^{{\bf m}_1+(2k+1){\bf m}_2}$ and
$t_0^{\bar 1}t^{{\bf m}_1+(2k+1){\bf m}_2}$ to  (3.14) respectively,
one gets that
\begin{eqnarray}\!\!\!\!\!\!\!\!\!\!\!\!&\!\!\!\!\!\!\!\!\!\!\!\!\!\!\!\!\!&
0=t_0^{\bar 0}t^{{\bf m}_1+(2k+1){\bf m}_2}\cdot
t_0^{\bar{0}}t^{-{\bf m}_1+2s{\bf m}_2}P_{e}(t^{{\bf m}_2})\cdot v_0
\nonumber\\[7pt]\!\!\!\!\!\!\!\!\!\!\!\!&\!\!\!\!\!\!\!\!\!\!\!\!\!\!\!\!\!\!\!\!\!\!\!\!&
\phantom{0}=q^{-m_{11}m_{12}-(2k+1)m_{11}m_{22}+2sm_{12}m_{21}+2s(2k+1)m_{21}m_{22}}\times
\nonumber\\[0pt]\!\!\!\!\!\!\!\!\!\!\!\!&\!\!\!\!\!\!\!\!\!\!\!\!\!\!\!\!\!\!\!\!\!\!\!\!&
\phantom{0=}\times\Big(t_0^{\bar 0}t^{(2k+2s+1){\bf
m}_2}P_{e}(t^{{\bf m}_2})-q^{(2k+2s+1)\alpha}t_0^{\bar
0}t^{(2k+2s+1){\bf
m}_2}P_{e}(q^{\alpha}t^{{\bf m}_2})\Big)\cdot v_0,%\mbox{ \ \ and}
\\[11pt]\!\!\!\!\!\!\!\!\!\!\!\!&\!\!\!\!\!\!\!\!\!\!\!\!\!\!\!\!\!\!\!\!\!\!\!\!&
0=t_0^{\bar 1}t^{{\bf m}_1+(2k+1){\bf m}_2}\cdot
(t_0^{\bar{0}}t^{-{\bf m}_1+2s{\bf m}_2}P_{e}(t^{{\bf m}_2}))\cdot
v_0
\nonumber\\[7pt]\!\!\!\!\!\!\!\!\!\!\!\!&\!\!\!\!\!\!\!\!\!\!\!\!\!\!\!\!\!\!\!\!\!\!\!\!&
\phantom{0}=q^{-m_{11}m_{12}-(2k+1)m_{11}m_{22}+2sm_{12}m_{21}+2s(2k+1)m_{21}m_{22}}\times
\nonumber\\[0pt]\!\!\!\!\!\!\!\!\!\!\!\!&\!\!\!\!\!\!\!\!\!\!\!\!\!\!\!\!\!\!\!\!\!\!\!\!&
\phantom{0=}\times\Big(t_0^{\bar 1}t^{(2k+2s+1){\bf
m}_2}P_{e}(t^{{\bf m}_2})-(-1)^{m_{11}}q^{(2k+2s+1)\alpha}t_0^{\bar
0}t^{(2k+2s+1){\bf m}_2}P_{e}(q^{\alpha}t^{{\bf m}_2})\Big)\cdot
v_0. \end{eqnarray} So we have
$$\begin{array}{lll}\Big(t_0^{\bar 0}t^{2k{\bf m}_2}P_{e}(t^{{\bf m}_2})-
q^{2k\alpha}t_0^{\bar 0}t^{2k{\bf m}_2}P_{e}(q^{\alpha}t^{{\bf
m}_2}) +a_{-k}q^{-4k^2m_{21}m_{22}}\beta\Big)\cdot v_0 =0,
\\[7pt]
\Big(t_0^{\bar 1}t^{2k{\bf m}_2}P_{e}(t^{{\bf m}_2})-(-1)^{m_{11}}
q^{2k\alpha}t_0^{\bar 1}t^{2k{\bf m}_2}P_{e}(q^{\alpha}t^{{\bf
m}_2})\Big)\cdot v_0 =0,
\\[7pt]
\Big(t_0^{\bar 0}t^{(2k+1){\bf m}_2}P_{e}(t^{{\bf
m}_2})-q^{(2k+1)\alpha}t_0^{\bar 0}t^{(2k+1){\bf
m}_2}P_{e}(q^{\alpha}t^{{\bf m}_2})\Big)\cdot v_0=0, \\[7pt]
\Big(t_0^{\bar 1}t^{(2k+1){\bf m}_2}P_{e}(t^{{\bf
m}_2})-(-1)^{m_{11}}q^{(2k+1)\alpha}t_0^{\bar 0}t^{(2k+1){\bf
m}_2}P_{e}(q^{\alpha}t^{{\bf m}_2})\Big)\cdot v_0=0, \end{array}$$
for all $k\in{\bf Z}$, which deduces the claim as required.
\vskip2mm

On the other hand,  we can choose an integer $s$ and a polynomial
$P_{o}(t^{{\bf m}_2})=\sum_{i=0}^{n}a_it^{2i{\bf m}_2}\in{\bf
C}[t^{{\bf m}_2}]$ with $a_0a_n\neq0$ such that
$$
t_0^{\bar{0}}t^{-{\bf m}_1+(2s+1){\bf m}_2}P_{o}(t^{{\bf m}_2})\cdot
v_0=0,
$$
since $\mbox{dim}V_{-1}<\infty$. Thus by a calculation similar to
the proof of Claim 1, we can deduce the following claim.\vskip2mm

{\bf Claim 2}\quad There is a polynomial $P_{o}(t^{{\bf
m}_2})=\sum_{i=0}^na_it^{2i{\bf m}_2}$ with $a_na_0\neq 0$ such that
\begin{eqnarray}&&\Big(t_0^{\bar 0}t^{2k{\bf m}_2}Q(t^{{\bf m}_2})-
q^{2k\alpha}t_0^{\bar 0}t^{2k{\bf m}_2}Q(q^{\alpha}t^{{\bf m}_2})
+a_{Q}\beta\Big)\cdot v_0 =0, \nonumber\\[0pt]&& \Big(t_0^{\bar 1}t^{2k{\bf m}_2}Q(t^{{\bf
m}_2})-(-1)^{m_{11}+1} q^{2k\alpha}t_0^{\bar 1}t^{2k{\bf
m}_2}Q(q^{\alpha}t^{{\bf m}_2})\Big)\cdot v_0 =0,
\nonumber\\[0pt]&&\Big(t_0^{\bar 0}t^{(2k+1){\bf m}_2}Q(t^{{\bf
m}_2})-q^{(2k+1)\alpha}t_0^{\bar 0}t^{(2k+1){\bf
m}_2}Q(q^{\alpha}t^{{\bf m}_2})\Big)\cdot v_0=0,
\nonumber\\[0pt]&&\Big(t_0^{\bar 1}t^{(2k+1){\bf m}_2}Q(t^{{\bf
m}_2})-(-1)^{m_{11}+1}q^{(2k+1)\alpha}t_0^{\bar 1}t^{(2k+1){\bf
m}_2}Q(q^{\alpha}t^{{\bf m}_2})\Big)\cdot v_0=0, \end{eqnarray}for
all $k\in{\bf Z}$ and $Q(t^{{\bf m}_2})\in{\bf C}[t^{\pm 2{\bf
m}_2}]$ with $P_{o}(t^{{\bf m}_2})\mid Q(t^{{\bf m}_2})$, where
$a_{Q}$ is the constant term of $t^{2k{\bf m}_2}Q(t^{{\bf
m}_2})$.\vskip2mm

Let $P(t^{{\bf m}_2})=\sum_{i=0}^na_it^{2i{\bf m}_2}$ be the product
of $P_o(t^{{\bf m}_2})$ and $P_e(t^{{\bf m}_2})$. We see that both
(3.13) and (3.19) hold for $P(t^{{\bf m}_2})$. Thus one can directly
deduce that both (3.10) and (3.12) hold for $P(t^{{\bf m}_2})$ and
$v_0\in V_0$. Since $v_0$ is a eigenvector of $t_0^{\bar 1}$, we
have
$$
0=t_0^{\bar 1}\cdot t_0^{\bar 1}t^{(2k+1){\bf m}_2}P(t^{{\bf
m}_2})\cdot v_0=[t_0^{\bar 1},t_0^{\bar 1}t^{(2k+1){\bf
m}_2}P(t^{{\bf m}_2})]\cdot v_0=2t_0^{\bar 0}t^{(2k+1){\bf
m}_2}P(t^{{\bf m}_2})\cdot v_0,
$$
and
$$
0=t_0^{\bar 1}\cdot t_0^{\bar 1}t^{(2k+1){\bf
m}_2}P(q^{\alpha}t^{{\bf m}_2})\cdot v_0=[t_0^{\bar 1},t_0^{\bar
1}t^{(2k+1){\bf m}_2}P(q^{\alpha}t^{{\bf m}_2})]\cdot v_0=2t_0^{\bar
0}t^{(2k+1){\bf m}_2}P(q^{\alpha}t^{{\bf m}_2})\cdot v_0,
$$
which deduces (3.11) for $P(t^{{\bf m}_2})$ and $v_0$.

From the definition of Lie subalgebra $L_0$, one can easily deduces
that if (3.10)--(3.12) hold for any $v\in V$, then they also hold
for $t_0^st^{k{\bf m}_2}\cdot v$, $\forall \ s\in{\bf Z}/2{\bf
Z},\;k\in{\bf Z}$. This completes the proof of this direction since
$V_0$ is an irreducible $L_0$-module.\vskip2mm

``$\Longleftarrow$''.\quad   We  first show the following claim by
induction on $s$.\vskip2mm

{\em {\bf Claim 3.}\quad For any $s\in{\bf Z}_+$ , there exists a
polynomial $P_s(t^{{\bf m}_2})=\sum_{j\in{\bf Z}}a_{s,j}t^{2j{\bf
m}_2}\in{\bf C}[t^{2{\bf m}_2}]$ such that
$$\begin{array}{lll}
\Big(t_0^{\bar 0}t^{2k{\bf m}_2}P_s(t^{{\bf
m}_2})-q^{2k\alpha}t_0^{\bar 0}t^{2k{\bf m}_2}P_s(q^{\alpha}t^{{\bf
m}_2})+a_{s,-k}q^{-4k^2m_{21}m_{22}}\beta\Big)\cdot V_{-s}=0,
\\[9pt]
t_0^{\bar 0}t^{(2k+1){\bf m}_2}P_s(t^{{\bf m}_2})\cdot
V_{-s}=t_0^{\bar 1}t^{k{\bf m}_2}P_s(t^{{\bf m}_2})\cdot V_{-s}=0,
\\[9pt]
t_0^rt^{-{\bf m}_1+k{\bf m}_2}P_s(t^{{\bf m}_2})\cdot
V_{-s}=0,\;\;\forall\ r\in{\bf Z}_2,k\in{\bf Z}. \end{array}$$ }

By the assumption and the definition of $L_0$-module $V_{0}$, one
can  deduce that the claim holds for $s=0$ with $P_0(t^{{\bf
m}_2})=P(t^{{\bf m}_2})$. Suppose it holds for $s$. Let us consider
the claim for $s+1$.

The equations in the claim are equivalent to
\begin{eqnarray}&&
\Big(t_0^{\bar 0}Q(t^{{\bf m}_2})-t_0^{\bar 0}Q(q^{\alpha}t^{{\bf
m}_2})+a_{Q}\beta\Big)\cdot V_{-s}=0, \nonumber\\[7pt]&&t_0^{\bar 0}t^{(2k+1){\bf
m}_2}Q(t^{{\bf m}_2})\cdot V_{-s}=t_0^{\bar 1}t^{k{\bf
m}_2}Q(t^{{\bf m}_2})\cdot V_{-s}=0,
\nonumber\\[7pt]&&t_0^rt^{-{\bf m}_1+k{\bf m}_2}Q(t^{{\bf m}_2})\cdot
V_{-s}=0,\;\;\forall\ r\in{\bf Z}_2,k\in{\bf Z}, \end{eqnarray}for
any $Q(t^{{\bf m}_2})\in{\bf C}[t^{\pm 2{\bf m}_2}]$ with
$P_s(t^{{\bf m}_2})\mid Q(t^{{\bf m}_2})$, where  $a_{Q}$ is the
constant term of $Q(t^{{\bf m}_2})$.

Let $P_{s+1}(t^{{\bf m}_2})=P_s(q^{\alpha}t^{{\bf m}_2})P_s(t^{{\bf
m}_2})P_s(q^{-\alpha}t^{{\bf m}_2})$. For any $p,r\in{\bf Z}_2,\,
j,k\in{\bf Z}$, using induction and by (3.20) we have
$$\begin{array}{lll}
\Big(t_0^{\bar 0}t^{2k{\bf m}_2}P_{s+1}(t^{{\bf
m}_2})-q^{2k\alpha}t_0^{\bar 0}t^{2k{\bf
m}_2}P_{s+1}(q^{\alpha}t^{{\bf
m}_2})+a_{s+1,-k}q^{-4k^2m_{21}m_{22}}\beta\Big)\cdot t_0^pt^{-{\bf
m}_1+j{\bf m}_2}\cdot V_{-s} \\[7pt]
=\Big[t_0^{\bar 0}t^{2k{\bf m}_2}P_{s+1}(t^{{\bf
m}_2})-q^{2k\alpha}t_0^{\bar 0}t^{2k{\bf
m}_2}P_{s+1}(q^{\alpha}t^{{\bf
m}_2})+a_{s+1,-k}q^{-k^2m_{21}m_{22}}\beta,t_0^pt^{-{\bf m}_1+j{\bf
m}_2}\Big]\cdot V_{-s}
\\[7pt]
=q^{2km_{22}(-m_{11}+jm_{21})}\Big(t_0^{p}t^{-{\bf m}_1+(2k+j){\bf
m}_2}\Big(P_{s+1}(q^{-\alpha}t^{{\bf
m}_2})-2q^{2k\alpha}P_{s+1}(t^{{\bf m}_2})
+q^{4k\alpha}P_{s+1}(q^{\alpha}t^{{\bf m}_2})\Big)\Big)\cdot
V_{-s}\\[7pt]
=0, \end{array}$$
 Thus, by (3.1), we obtain that
$$
\Big(t_0^{\bar 0}t^{2k{\bf m}_2}P_{s+1}(t^{{\bf
m}_2})-q^{2k\alpha}t_0^{\bar 0}t^{2k{\bf
m}_2}P_{s+1}(q^{\alpha}t^{{\bf
m}_2})+a_{s+1,-k}q^{-4k^2m_{21}m_{22}}\beta\Big)\cdot V_{-s-1}=0.
\eqno{(3.21)}
$$
Similarly, one can prove that
$$t_0^{\bar 0}t^{(2k+1){\bf m}_2}P_{s+1}(t^{{\bf m}_2})\cdot V_{-s-1}=t_0^{\bar
1}t^{k{\bf m}_2}P_{s+1}(t^{{\bf m}_2})\cdot V_{-s-1}=0,\;\; \forall\
k\in{\bf Z}. \eqno{(3.22)}
$$
 This proves the first two equations in the
claim for $s+1$.

Using (3.21), (3.22)  and induction, we deduce that for any
$l,k\in{\bf Z},\,\,n,r\in{\bf Z}_2$,
$$\begin{array}{lll}
t_0^nt^{{\bf m}_1+l{\bf m}_2}\cdot t_0^rt^{-{\bf m}_1+k{\bf
m}_2}P_{s+1}(t^{{\bf m}_2})\cdot V_{-s-1} \\[7pt]=[t_0^nt^{{\bf
m}_1+l{\bf m}_2},t_0^rt^{-{\bf m}_1+k{\bf m}_2}P_{s+1}(t^{{\bf
m}_2})]\cdot V_{-s-1}+ t_0^rt^{-{\bf m}_1+k{\bf m}_2}P_{s+1}(t^{{\bf
m}_2})\cdot t_0^nt^{{\bf m}_1+l{\bf m}_2}\cdot V_{-s-1}
\\[7pt]=(-1)^{r(m_{11}+lm_{21})}q^{-m_{11}m_{12}+km_{12}m_{21}-lm_{11}m_{22}+lkm_{21}m_{22}}\Big(t_0^{n+r}t^{(l+k){\bf
m}_2}P_{s+1}(t^{{\bf m}_2})
\\[7pt]\phantom{=}-(-1)^{(n+r)m_{11}+nk+rl}q^{(k+l)\alpha}t_0^{n+r}t^{(l+k){\bf
m}_2}P_{s+1}(q^{\alpha}t^{{\bf
m}_2})\\[7pt]\phantom{=}+a_{s+1,i}\delta_{k+l+2i,0}\delta_{r+n,{\bar
0}}q^{-(l+k)^2m_{21}m_{22}}\beta\Big)\cdot V_{-s-1}\\[7pt]=0. \end{array}$$
Hence, by (3.2),
$$t_0^rt^{-{\bf m}_1+k{\bf
m}_2}P_{s+1}(t^{{\bf m}_2})\cdot  V_{-s-1}=0,
$$
 for all $r\in{\bf
Z}_2,\;k\in{\bf Z}$, which implies the third equation in the claim
for $s+1$. Therefore the claim follows by induction.\vskip2mm

From the third equation of the claim and (3.1), we see that
$$
\mbox{dim}V_{-s-1}\leq 2\mbox{deg}(P_{s+1}(t^{{\bf m}_2}))\cdot
\mbox{dim}V_{s},\;\;\forall \ s\in{\bf Z}_+,
$$
where $\mbox{deg}(P_{s+1}(t^{{\bf m}_2}))$ denotes the degree of
polynomial $P_{s+1}(t^{{\bf m}_2})$. Hence $\OVERLINE{M}^+(V(\ula
,\psi),{\bf m}_1,{\bf m}_2)$ $\in{\cal O}_{\bf Z}$.
 \hfill$\Box$\vskip2mm

{\em{\bf Theorem 3.4}\quad Let $m_{21}$ be an odd integer. Then
$\OVERLINE{M}^+(\ula ,\psi,{\bf m}_1,{\bf m}_2)\in{\cal O}_{\bf Z}$
if and only if there exist
$b_{10},b_{11},\cdots,b_{1s_1},b_{20},b_{21},\cdots,b_{2s_2},\cdots,b_{r0},b_{r1},\cdots,b_{rs_r}\in{\bf
C}\,$ and $\alpha_1,\cdots,\alpha_r\in{\bf C}^{*}$ such that for any
$i\in{\bf Z}^{*}$, $j\in{\bf Z}_2$,
$$\begin{array}{ll}\displaystyle
\psi(t_0t^{2i{\bf
m}_2})=\frac{(b_{10}+b_{11}i+\cdots+b_{1s_1}i^{s_1})\alpha_1^i
+\cdots+(b_{r0}+b_{r1}i+\cdots,b_{rs_r}i^{s_r})\alpha_r^i}{(1-q^{2i\alpha
})q^{2i^2m_{21}m_{22}}}, \\[11pt]\psi(\beta)=b_{10}+b_{20}+\cdots+b_{r0},
\mbox{ \ \ \ \ and \ \ \ \
}\psi(m_{21}c_1+m_{22}c_2)=0,\end{array}$$
 where $\alpha=m_{11}m_{22}-m_{21}m_{12}\in\{\pm
1\}$.}\vskip2mm

{\bf Proof } ``$\Longrightarrow$''. Let
$f_i=\psi((1-q^{2i\alpha})q^{2i^2m_{21}m_{22}}t_0^{\bar 0}t^{2i{\bf
m}_2})$ for $i\in{\bf Z}^*$ and $f_0=\psi(\beta)$. By Lemma 3.3,
there exist complex numbers $a_0,a_1,\cdots,a_n$ with $a_0a_n\neq 0$
such that
$$
\sum\limits_{i=0}^na_iq^{-2i^2m_{21}m_{22}}f_{k+i}=0,\;\forall\
k\in{\bf Z}.
$$
Thus, by using a technique in the proof of Theorem 3.2, we can
deduce the result as required.\vskip2mm

``$\Longleftarrow$''. Set
$$
Q(x)=\Big(\prod\limits_{i=1}^r(x-\alpha_i)^{s_i+1}\Big)\Big(\prod\limits_{j=1}^{\nu}(x-a_j)\Big)\Big
(\prod\limits_{j=1}^{\nu}(x-q^{2\alpha}a_j)\Big)=:\sum\limits_{i=1}^nb_ix^i,
$$
and
$$
 f_i=\psi\Big((1-q^{2i\alpha})q^{2i^2m_{21}m_{22}}t_0^{\bar
0}t^{2i{\bf m}_2}\Big),\;\forall\ i\in{\bf Z}^*,\quad
f_0=\psi(\beta).
$$
Then one can easily verify that
$$
\sum\limits_{i=0}^nb_if_{k+i}=0,\quad \forall\ k\in{\bf Z}.
\eqno{(3.23)}
$$
Meanwhile, we have $(\prod_{j=1}^{\nu}(x-a_j))\mid x^kQ(x)$ and
$(\prod_{j=1}^{\nu}(x-a_j))\mid x^kQ(q^{2\alpha}x)$ for any
$k\in{\bf Z}$, which deduces\setcounter{equation}{23}
\begin{eqnarray}&&
\sum\limits_{i=1}^nb_iq^{\frac{1}{2}(2i+2k+1)^2m_{22}m_{21}}t_0^st^{(2i+2k+1){\bf
m}_2}\cdot V_0=0,\\&&
\sum\limits_{i=1}^nb_iq^{2i\alpha}q^{\frac{1}{2}(2i+2k+1)^2m_{22}m_{21}}t_0^st^{(2i+2k+1){\bf
m}_2}\cdot V_0=0, \quad \forall\ s\in{\bf Z}_2,
\end{eqnarray}and\begin{eqnarray}%\\
&& \sum\limits_{i=1}^nb_iq^{2(i+k)^2m_{22}m_{21}}t_0^{\bar
1}t^{2(i+k){\bf m}_2}\cdot V_0=0,
\\&& \sum\limits_{i=1}^nb_iq^{2i\alpha}q^{2(i+k)^2m_{22}m_{21}}t_0^{\bar
1}t^{2(i+k){\bf m}_2}\cdot V_0=0, \end{eqnarray} by Remark 2.6. Let
$b_i'=q^{2i^2m_{21}m_{22}}b_i$ for $0\leq i\leq n$ and
$P(x)=\sum_{i=1}^nb_i'x^i$. By (3.23) and the construction of
$V(\ula ,\psi)$, we have
$$\begin{array}{lll}
\Big(t_0^{\bar 0}t^{2k{\bf m}_2}P(t^{2{\bf
m}_2})-q^{2k\alpha}t_0^{\bar 0}t^{2k{\bf m}_2}P(q^{2\alpha}t^{{2\bf
m}_2})+b'_{-k}q^{-4k^2m_{21}m_{22}}\beta\Big)\cdot V_0
\\[7pt]=q^{-2k^2m_{21}m_{22}}\psi\Big(\sum\limits_{i=1}^nb_i(1-q^{2(k+i)\alpha})q^{2(k+i)^2m_{22}m_{21}}t_0^{\bar
0}t^{2(k+i){\bf m}_2}+b_{-k}\beta\Big)\cdot V_0
\\[7pt]=q^{-2k^2m_{21}m_{22}}\sum\limits_{i=1}^nb_if_{k+i}\cdot V_0\\[7pt]=0,
\end{array}$$ which deduces (3.10). Similarly, we have, for any $k\in{\bf
Z}$,
$$\begin{array}{ll}
t_0^{s}t^{(2k+1){\bf m}_2}P(t^{2{\bf m}_2})\cdot
V_0\!\!\!\!&=\sum\limits_{i=1}^{n}b_iq^{(2i^2+4ki+2i)m_{21}m_{22}}t_0^st^{(2k+2i+1){\bf
m}_2}\cdot V_0
\\[7pt]&=q^{-2k^2-2k-\frac{1}{2}}\sum\limits_{i=1}^nb_iq^{\frac{1}{2}(2k+2i+1)^2m_{21}m_{22}}t_0^st^{(2k+2i+1){\bf
m}_2}\cdot V_0\\[7pt]&=0,
\end{array}$$ and
$$\begin{array}{ll}
t_0^{s}t^{(2k+1){\bf m}_2}P(q^{2\alpha}t^{2{\bf m}_2})\cdot
V_0\!\!\!\!&=\sum\limits_{i=1}^{n}b_iq^{2i\alpha+(2i^2+4ki+2i)m_{21}m_{22}}t_0^st^{(2k+2i+1){\bf
m}_2}\cdot V_0
\\[7pt]&=q^{-2k^2-2k-\frac{1}{2}}\sum\limits_{i=1}^nb_iq^{2i\alpha}q^{\frac{1}{2}(2k+2i+1)^2m_{21}m_{22}}t_0^st^{(2k+2i+1){\bf
m}_2}\cdot V_0\\[7pt]&=0,
\end{array}$$ by (3.24) and (3.25) respectively. Now one can easily deduce
the following equation
$$
t_0^{\bar 1}t^{2k{\bf m}_2}P(t^{2{\bf m}_2})\cdot V_0=0, \mbox{ \ \
\ and \ \ \ }t_0^{\bar 1}t^{2k{\bf m}_2}P(q^{2\alpha}t^{2{\bf
m}_2})\cdot V_0=0,
$$
by using (3.26) and (3.27) respectively. Therefore (3.10)--(3.12)
hold for $P(t^{2{\bf m}_2})=\sum_{i=1}^nb_i't^{2i{\bf m}_2}$. Thus
$\OVERLINE{M}^+(\ula ,\psi,{\bf m}_1,{\bf m}_2)\in{\cal O}_{\bf Z}$
by Lemma 3.3. \hfill$\Box$\vskip2mm

{\em {\bf Remark 3.5 } A linear function $\psi$ over $L_0$ having
the form as described in Theorem 3.2 is called  an exp-polynomial
function over $L_0$. Similarly, a linear function $\psi$ over ${\cal
A}$ having the form as described in Theorem 3.4 is called an
exp-polynomial function over ${\cal A}$.} \vspace{7pt}

\centerline{\large\bf \S 4 \ \ Classification of the generalized
highest }\centerline{\large\bf weight irreducible {\bf Z}-graded
$L$-modules}

 { \em{\bf Lemma
4.1}\quad If $V$ is a nontrivial irreducible generalized highest
weight {\bf Z}-graded $L$-module corresponding to a ${\bf Z}$-basis
$B=\{{\bf b}_1, {\bf b}_2\}$ of ${\bf Z}^2$, \vspace*{-4pt}then
\begin{itemize}\parskip-5pt\item[$(1)$]  For any $v\in V$ there is some $p\in{\bf N}$ such that
$t_0^it^{m_1{\bf b}_1+m_2{\bf b}_2}\cdot  v=0$ for all $m_1,m_2\geq
p$ and $i\in{\bf Z}_2$.
\item[$(2)$]  For any $0\neq v\in V$ and $m_1,m_2>0,\,i\in{\bf Z}_2$, we
have $t_0^it^{-m_1{\bf b}_1-m_2{\bf b}_2}\cdot  v\neq
0$.\end{itemize}}\vskip1mm

{\bf Proof }\quad Assume that $v_0$ is a generalized highest weight
vector corresponding to the ${\bf Z}$-basis $B=\{{\bf b}_1, {\bf
b}_2\}$ of ${\bf Z}^2$.

(1) By the irreducibility of $V$ and the PBW theorem, there exists
$u\in U(L)$ such that $v=u\cdot  v_0$, where  $u$ is a linear
combination of elements of the form $$u_n=(t_0^{k_1}t^{i_1{\bf
b}_1+j_1{\bf b}_2})\cdot(t_0^{k_2}t^{i_2{\bf b}_1+j_12{\bf
b}_2})\cdots (t_0^{k_n}t^{i_n{\bf b}_1+j_n{\bf b}_2}),$$ where,
``\,$\cdot$\,'' denotes the product in $U(L)$. Thus, we may assume
$u=u_n$. Take
$$
p_1=-\sum_{i_s<0}i_s+1,\quad p_2=-\sum_{j_s<0}j_s+1.
$$
By induction on $n$, one gets that $t_0^kt^{i{\bf b}_1+j{\bf
b}_2}\cdot v=0$ for any $k\in{\bf Z}_2,i\geq p_1$ and $j\geq p_2$,
which gives the result with $p=\mbox{max}\{p_1,p_2\}$.

(2) Suppose there are $0\neq v\in V$ and $i\in{\bf Z}_2,m_1,m_2>0$
with
$$t_0^it^{-m_1{\bf b}_1-m_2{\bf b}_2}\cdot  v=0.
$$
Let $p$ be as in the proof of (1). Then
$$t_0^it^{-m_1{\bf
b}_1-m_2{\bf b}_2},\;t_0^jt^{{\bf b}_1+p(m_1{\bf b}_1 +m_2{\bf
b}_2)},\;t_0^jt^{{\bf b}_2+p(m_1{\bf b}_1+m_2{\bf b}_2)},
\;\;\forall j\in{\bf Z}_2,
$$
act trivially on $v$. Since the above elements generate the Lie
algebra $L$. So $V$ is a trivial module, a
contradiction.\hfill$\Box$

\vspace{3mm}{ \em{\bf Lemma 4.2}\quad If $V\in{\cal O}_{\bf Z}$ is a
generalized highest weight $L$-module corresponding to the ${\bf
Z}$-basis $B=\{{\bf b}_1, {\bf b}_2\}$ of ${\bf Z}^2$, then $V$ must
be a highest or  lowest weight module.}\vskip2mm

{\bf Proof }\quad Suppose $V$ is a generalized highest weight module
corresponding to the ${\bf Z}$-basis $\{{\bf b}_1=b_{11}{\bf
m}_1+b_{12}{\bf m}_2$, ${\bf b}_2=b_{21}{\bf m}_1+b_{22}{\bf m}_2\}$
of ${\bf Z}^2$. By shifting index of $V_i$ if necessary, we can
suppose the highest degree of $V$ is $0$. Let $a=b_{11}+b_{21}$ and
$$
\wp(V)=\{m\in {\bf Z}\mid V_m\neq 0\}.
$$
We may assume $a\neq 0$: In fact, if $a=0$ we can choose ${\bf
b}_1'=3{\bf b}_1+{\bf b}_2,\ {\bf b}_2'=2{\bf b}_1+{\bf b}_2$, then
$V$ is a generalized highest weight {\bf Z}-graded module
corresponding to the ${\bf Z}$-basis $\{{\bf b}_1',{\bf b}_2'\}$ of
${\bf Z}^2$. Replacing ${\bf b}_1,{\bf b}_2$ by ${\bf b}_1',{\bf
b}_2'$ gives $a\neq 0$.

Now we prove that if $a>0$ then $V$ is a highest weight module.
 Let
$$
{\cal A}_i=\{j\in{\bf Z}\mid i+aj\in\wp(V)\},\;\; \forall\  0\leq
i<a.
$$
 Then there is  $m_i\in{\bf Z}$ such that ${\cal
A}_i=\{j\in{\bf Z}\mid j\leq m_i\} $ or ${\cal A}_i={\bf Z}$ by
Lemma 4.1(2).

Set ${\bf b}={\bf b}_1+{\bf b}_2$. We want to prove ${\cal
A}_i\not={\bf Z}$ for all $0\leq i<a$. Otherwise, (by shifting the
index of ${\cal A}_i$ if necesary) we may assume ${\cal A}_0={\bf
Z}$. Thus we can choose $0\neq v_j\in V_{aj}$ for any $j\in{\bf Z}$.
By Lemma 4.1(1), we know that there is
 $p_{v_j}>0$ with
$$
t_0^kt^{s_1{\bf b}_1+s_2{\bf b}_2}\cdot v_j=0,\;\;\forall\
s_1,s_2>p_{v_j},\; k\in{\bf Z}/2{\bf Z}. \eqno{(4.1)}
$$
Choose $\{k_j\in{\bf N}\mid j\in{\bf N}\}$ and $v_{k_j}\in V_{ak_j}$
such that
$$
k_{j+1}>k_j+p_{v_{k_j}}+2. \eqno{(4.2)}
$$
We prove that $\{t_0^{{\bar 0}}t^{-k_j{\bf b}}\cdot v_{k_j}\mid
j\in{\bf N}\}\subset V_0$ is a set of linearly independent vectors,
from which we can get a contradiction and thus deduces the result as
we hope.

Indeed, for any $r\in{\bf N}$, there exists $a_r\in{\bf N}$ such
that $t_0^0t^{x{\bf b}+{\bf b}_1}v_{k_r}=0,\;\forall x\geq a_r$ by
Lemma 4.1(1). On the other hand, we know that $t_0^0t^{x{\bf b}+{\bf
b}_1}\cdot v_{k_r}\neq 0$ for any $x<-1$ by Lemma 4.1(2). Thus we
can choose $s_r\geq -2$ such that
$$
t_0^{{\bar 0}}t^{s_r{\bf b}+{\bf b}_1}\cdot v_{k_r}\not=0,\quad\
\quad t_0^{{\bar 0}}t^{x{\bf b}+{\bf b}_1}\cdot v_{k_r}=0,\;\forall
x>s_r. \eqno{(4.3)}
$$
By (4.2) we have $k_r+s_r-k_j>p_{v_j}$ for all $1\leq j<r$. Hence by
(4.1) we know that for all $1\leq j<r$,
$$\begin{array}{ll}
t_0^{{\bar 0}}t^{(k_r+s_r){\bf b}+{\bf b}_1}\cdot t_0^{{\bar
0}}t^{-k_j{\bf b}}\cdot v_{k_j}\!\!\!\!\\[7pt] =[t_0^{{\bar
0}}t^{(k_r+s_r){\bf b}+{\bf b}_1},t_0^{{\bar 0}}t^{-k_j{\bf
b}}]\cdot v_{k_j}
\\[7pt]=q^{-k_j((k_r+s_r)(b'_{12}+b'_{22})+b'_{12})(b'_{11}+b'_{21})}(1-q^{k_j(b'_{12}b'_{21}-b'_{11}b'_{22})})
t_0^{{\bar 0}}t^{(k_r+s_r-k_j) {\bf b}+{\bf b}_1}\cdot v_{k_j}\\[7pt]=0,
\end{array}$$ where $$b'_{11}=b_{11}m_{11}+b_{12}m_{21},\
b'_{12}=b_{11}m_{12}+b_{12}m_{22},\
b'_{21}=b_{21}m_{11}+b_{22}m_{21},\
b'_{22}=b_{21}m_{12}+b_{22}m_{22}.$$  Now by (4.2) and (4.3), one
gets
$$\begin{array}{ll}
t_0^{{\bar 0}}t^{(k_r+s_r){\bf b}+{\bf b}_1}\cdot t_0^{{\bar
0}}t^{-k_r{\bf b}}\cdot v_{k_r}\\[7pt]
=[t_0^{{\bar 0}}t^{(k_r+s_r){\bf b}+{\bf b}_1},t_0^{{\bar
0}}t^{-k_r{\bf b}}]\cdot v_{k_r}
\\[7pt]
=q^{-k_r((k_r+s_r)(b'_{12}+b'_{22})+b'_{12})(b'_{11}+b'_{21})}(1-q^{k_r(b'_{12}b'_{21}-b'_{11}b'_{22})})t_0^{{\bar
0}}t^{s_r {\bf b}+{\bf b}_1}\cdot v_{k_r}\\[7pt]
\not=0. \end{array}$$ Hence if $\sum_{j=1}^{n}\lambda_jt_0^{{\bar
0}}t^{-k_j{\bf b}}\cdot v_{k_j}=0$ then $\lambda_n=\lambda_{n-1}
=\cdots=\lambda_{1}=0$ by the arbitrariness of $r$. So we see that
$\{t_0^{{\bar 0}}t^{-k_j{\bf b}}\cdot v_{k_j}\mid j\in{\bf
N}\}\subset V_0$ is a set of linearly independent vectors, which
contradicts the fact that $V\in{\cal O}_{\bf Z}$. Therefore, for any
$0\leq i<a$, there is $m_i\in{\bf Z}$ such that
 ${\cal A}_i=\{j\in{\bf Z}\mid j\leq m_i\} $, which deduces that
 $V$ is a highest weight module since
$\wp(V)=\bigcup_{i=0}^{a-1}{\cal A}_i$.

Similarly, one can prove that if $a<0$ then $V$ is a lowest weight
module.\hfill$\Box$\vskip2mm

From the above lemma and the results in Section 3, we have the
following theorem.\vskip2mm

{\em {\bf Theorem 4.3}  $V$ is a quasi-finite irreducible ${\bf
Z}$-graded $L$-module if and only if one of the following statements
\vspace*{-4pt}hold:
\begin{itemize}\parskip-5pt\item[$(1)$]
$V$ is a uniformly bounded module;
\item[$(2)$] If $m_{21}$ is an even integer then there exists an
exp-polynomial function $\psi$ over $L_0$ such that $$V\cong
\OVERLINE{M}^+(\psi,{\bf m}_1,{\bf m}_2)\mbox{ \ \  or \ \ }V\cong
\OVERLINE{M}^{ -}(\psi,{\bf m}_1,{\bf m}_2);$$
\item[$(3)$] If $m_{21}$ is an odd integer then there exist an exp-polynomial
function $\psi$ over ${\cal A}$, a finite sequence of nonzero
distinct numbers $\ula =(a_1,\cdots,a_{\nu})$ and some finite
dimensional irreducible $sl_2$-modules $V_1,\cdots,V_{\nu}$ such
that $$V\cong \OVERLINE{M}^+(\ula ,\psi,{\bf m}_1,{\bf m}_2)\mbox{ \
\ or  \ \ }V\cong \OVERLINE{M}^{ -}(\ula ,\psi,{\bf m}_1,{\bf
m}_2).$$\end{itemize}}
 \vspace{1mm}{ \it{\bf Theorem 4.4~~(Main Theorem)}\quad If $V$ is a quasi-finite irreducible
 {\bf Z}-graded $L$-module with nontrivial center
 then one of the following statements must \vspace*{-4pt}hold:
\begin{itemize}\parskip-5pt\item[$(1)$]
If $m_{21}$ is an even integer then there exists an exp-polynomial
function $\psi$ over $L_0$ such that
$$V\cong
\OVERLINE{M}^+(\psi,{\bf m}_1,{\bf m}_2)\mbox{ \ \  or \ \ }V\cong
\OVERLINE{M}^{ -}(\psi,{\bf m}_1,{\bf m}_2);$$
\item[$(2)$] If $m_{21}$ is an odd integer then there exist an exp-polynomial
function $\psi$ over ${\cal A}$, a finite sequence of nonzero
distinct numbers $\ula =(a_1,\cdots,a_{\nu})$ and some finite
dimensional irreducible $sl_2$ modules $V_1,\cdots,V_{\nu}$ such
that
$$V\cong \OVERLINE{M}^+(\ula ,\psi,{\bf m}_1,{\bf
m}_2)\mbox{ \ \ or  \ \ }V\cong \OVERLINE{M}^{ -}(\ula ,\psi,{\bf
m}_1,{\bf m}_2).$$\end{itemize}}

 {\bf Proof£º}\quad By Theorem 4.3, we only need to show that $V$ is not a uniformly bounded module.
 From the definition of Lie algebra $L$, we see that ${\cal H}_i=\la t_0^{\bar 0}t^{k{\bf m}_i}, m_{i1}c_1+m_{i2}c_2
 \mid k\in{\bf Z}^*\ra ,\; i=1,2$ are Heisenberg Lie algebras.
 As $V$ is a quasi-finite irreducible {\bf Z}-graded $L$-module, we deduce that $m_{21}c_1+m_{22}c_2$
 must be zero. Thus, by the assumption, we have that
 $m_{11}c_1+m_{12}c_2\neq 0$ since
 $\{{\bf m}_1,{\bf m}_2\}$ is a ${\bf Z}$-basis of ${\bf Z}^2$.
 Therefore, $V$ is not a uniformly bounded module by a well-known result about the representation
 of the Heisenberg Lie algebra. \hfill$\Box$\vskip2mm

We  close this section by showing that nontrivial modules
$\OVERLINE{M}^+(\psi,{\bf m}_1,{\bf m}_2),\,M^+(\ula ,\psi,{\bf
m}_1,{\bf m}_2) $ are not uniformly bounded and not
integrable.\vskip2mm

{ \em {\bf Theorem 4.5}\quad Nontrivial module
$\OVERLINE{M}^+(\psi,{\bf m}_1,{\bf m}_2)$ or $M^+(\ula ,\psi,{\bf
m}_1,{\bf m}_2) $ is not uniformly bounded.}\vskip2mm

 {\bf Proof£º}\quad Set $V\cong \OVERLINE{M}^+(\psi,{\bf
m}_1,{\bf m}_2)$ or $V\cong \OVERLINE{M}^+(\ula ,\psi,{\bf m}_1,{\bf
m}_2)$ and $V=\oplus_{k\in{\bf Z}_+}V_{-k}$. Since $V$ is not
trivial, there exist $v_0\in V_0$, $k\in {\bf Z}^*$ and $l\in {\bf
Z}_2$ such that $t_0^lt^{k{\bf m}_2}\cdot v_0\neq 0$. Thus
$$\begin{array}{ll}
t_0^{\bar 0}t^{{\bf m}_1}\cdot t_0^lt^{-{\bf m}_1+k{\bf m}_2}\cdot
v_0\!\!\!\!&=[t_0^{\bar 0}t^{{\bf m}_1},t_0^lt^{-{\bf m}_1+k{\bf
m}_2}]v_0
\\[7pt]&=((-1)^{lm_{11}}q^{m_{12}(-m_{11}+km_{21})}-q^{m_{11}(-m_{12}+km_{22})})t_0^lt^{k{\bf
m}_2}\cdot v_0\neq 0, \end{array}$$ which deduces that
$t_0^lt^{-{\bf m}_1+k{\bf m}_2}\cdot v_0\neq 0$.

Next, we prove that if $0\neq v_{-m}\in V_{-m}$ then $t_0^{\bar
0}t^{-{\bf m}_1}\cdot v_{-m}\neq 0$. Suppose $t_0^{\bar 0}t^{-{\bf
m}_1}\cdot v_{-m}=0$ for some $0\neq v_{-m}\in V_{-m}$. From the
construction of $V$, we know that $t_0^lt^{(m+1){\bf m}_1\pm {\bf
m}_2}$ also act trivially on $v_{-m}$ for any $l\in{\bf Z}_2$. Since
$L$ is generated by $ t_0^{\bar 0}t^{-{\bf m}_1},\,
t_0^lt^{(m+1){\bf m}_1\pm{\bf m}_2},\, l=\bar 0,\bar 1$, we see that
$V$ is a trivial module, a contradiction.

Set
$$
{\mathscr A}=\{(t_0^{\bar 0}t^{-{\bf m}_1})^j\cdot
t_0^{l}t^{(-n+j){\bf m}_1+k{\bf m}_2}\cdot v_0\mid 0\leq j<
n\}\subset V_{-n}, \forall\ n\in{\bf N}.
 $$
Now we prove that ${\mathscr A}$ is a set of linear independent
vectors. If
 $$\sum_{j=0}^{n-1}\lambda_j(t_0^{\bar 0}t^{-{\bf m}_1})^jt_0^{l}t^{(-n+j){\bf m}_1+k{\bf m}_2}\cdot v_0=0,
 $$
 then for any $0\leq i< n-1$ we have
 $$\begin{array}{ll}
 0\!\!\!\!&
 =q^{n(n-i)m_{11}m_{12}-k(n-i)m_{12}m_{21}}t_0^{\bar 0}t^{(n-i){\bf m}_1}\cdot \sum\limits_{j=0}^{n-1}\lambda_j
 (t_0^{\bar 0}t^{-{\bf m}_1})^j\cdot t_0^{l}t^{(-n+j){\bf m}_1+k{\bf m}_2}\cdot v_0
\\[7pt]&
  =\sum\limits_{j=0}^{i}\lambda_jq^{j(n-i)m_{11}m_{12}}
  ((-1)^{l(n-i)m_{11}}-q^{k(n-i)\alpha})(t_0^{\bar 0}t^{-{\bf m}_1})^j\cdot t_0^{l}t^{(j-i){\bf m}_1+k{\bf
  m}_2}\cdot v_0,
\end{array} $$
 where $\alpha=m_{11}m_{22}-m_{12}m_{21}$, which deduces $\lambda_0=\cdots=\lambda_{n-1}=0$.
 Hence ${\mathscr A}$ is a set of linear independent vectors
 in $V_{-n}$ and thus
 $$
 \mbox{dim}V_{-n}\geq n.
 $$
Therefore $V$ is not a uniformly bounded module by the arbitrariness
of $n$. \hfill$\Box$

\vspace{2mm} In [21], Rao gives a classification of the integrable
modules with nonzero center for the core of EALAs coordinatized by
quantum tori. We want to prove that the $L$-modules constructed in
this paper are in general not integrable. First we recall the
concept of the integrable modules. Let $\tau$ be the Lie algebra
defined in Section 2.  A $\tau$-module $V$ is {\it integrable} if,
for any $v\in V$ and ${\bf m}\in{\bf Z}^2$, there exist
$k_1=k_1({\bf m},v),k_2=k_2({\bf m},v)$ such that $(E_{12}(t^{\bf
m}))^{k_1}\cdot v=(E_{21}(t^{\bf m}))^{k_2}\cdot v=0$. Thus by
Proposition 2.1, an $L$-module $V$ is integrable if, for any $v\in
V$ and ${\bf m}=(2m_1+1,m_2)\in{\bf Z}^2$, there exist $k_1=k_1({\bf
m},v),k_2=k_2({\bf m},v)$ such that $$ (t_0^{\bar 0}t^{\bf
m}+t_0^{\bar 1}t^{\bf m})^{k_1}\cdot v=0=(t_0^{\bar 0}t^{\bf
m}-t_0^{\bar 1}t^{\bf m})^{k_2}\cdot v=0.$$

\vspace{2mm} { \em {\bf Theorem 4.6}\quad Nontrivial modules
$\OVERLINE{M}^+(\psi,{\bf m}_1,{\bf m}_2)$ or $\OVERLINE{M}^+(\ula
,\psi,{\bf m}_1,{\bf m}_2)$ is not integrable.}\vskip2mm

{\bf Proof} Set $V\cong \OVERLINE{M}^+(\psi,{\bf m}_1,{\bf m}_2)$ or
$V\cong \OVERLINE{M}^+(\ula ,\psi,{\bf m}_1,{\bf m}_2)$ and
$V=\oplus_{k\in{\bf Z}_+}V_{-k}$. Choose two positive integers $a$
and $b$ such that ${\bf m}=a{\bf m}_1+b{\bf m}_2=:(2k+1,l)$. Let
$v_0\in V_0$ be an eigenvector of $t_0^{\bar 1}$. Then we have
$$
(t_0^{\bar 0}t^{\bf m}\pm t_0^{\bar 1}t^{\bf m})\cdot v_0=0,
$$
by the construction of $V$. On the other hand, by using the
isomorphism $\varphi$ defined in Proposition 2.1, we have
$$
\varphi(t_0^{\bar 0}t^{\bf m}+ t_0^{\bar 1}t^{\bf
m})=2E_{21}(t_1^{m_1+1}t_2^{m_2}),\quad \varphi(t_0^{\bar 0}t^{\bf
m}- t_0^{\bar 1}t^{\bf m})=2q^{-m_2}E_{12}(t_1^{m_1}t_2^{m_2}),
$$
and
$$
\varphi(t_0^{\bar 0}t^{-\bf m}+ t_0^{\bar 1}t^{-\bf
m})=2E_{21}(t_1^{-m_1}t_2^{-m_2}),\quad \varphi(t_0^{\bar 0}t^{-\bf
m}- t_0^{\bar 1}t^{-\bf
m})=2q^{m_2+1}E_{12}(t_1^{-m_1-1}t_2^{-m_2}).
$$
Thus, by a well-known result on the $sl_2$-modules, we see that if
$V$ is integrable then
$$
t_0^{\bar 1}\cdot v_0=0, \quad (t_0^{\bar 0}t^{-\bf m}+ t_0^{\bar
1}t^{-\bf m})\cdot v_0=0, \quad (t_0^{\bar 0}t^{-\bf m}- t_0^{\bar
1}t^{-\bf m})\cdot v_0=0.
$$
So  $t_0^{\bar 0}t^{-\bf m},\,t_0^{\bar 1}t^{-\bf m}$ act trivially
on $v_0$. On the other hand, the construction of $V$ shows that
$t_0^it^{2{\bf m}\pm{\bf m}_1},\,t_0^it^{2{\bf m}\pm{\bf m}_2}$ act
trivially on $v_0$. Thus $L$ acts trivially on $v_0$ since $L$ is
generated by $t_0^{\bar 0}t^{-\bf m},\, t_0^{\bar 1}t^{-\bf
m},\,t_0^it^{2{\bf m}\pm{\bf m}_1},\,t_0^it^{2{\bf m}\pm{\bf m}_2}$.
Hence $V$ is a trivial $L$-module, a contradiction.
\hfill$\Box$\vskip7pt

\centerline{\large\bf \S 5 \ \ Two classes of highest weight ${\bf
Z}^2$-graded $L$-modules}

In this section, we  construct two classes of highest weight
quasi-finite irreducible ${\bf Z}^2$-graded $L$-modules. For any
highest weight {\bf Z}-graded $L$-module $V=\oplus_{k\in{\bf
Z}_+}V_{-k}$, set $V_{{\bf Z}^2}=V\otimes{\bf C}[x^{\pm 1}]$. We
define the action of the elements of $L$ on $V_{{\bf Z}^2}$ as
follows
$$
t_0^it^{m{\bf m}_1+n{\bf m}_2}\cdot (v\otimes x^r)=(t_0^it^{m{\bf
m}_1+n{\bf m}_2}\cdot v)\otimes x^{r+n},
$$
where $v\in V,\ i\in{\bf Z}_2,\ r,m,n\in{\bf Z}$. For any $v_{-k}\in
V_{-k}$, we define the degree of $v_{-k}\otimes t^r$ to be $-k{\bf
m}_1+r{\bf m}_2$. Then one can easily see that $V_{{\bf Z}^2}$
becomes a ${\bf Z}^2$-graded $L$-module. Let $W$ be an irreducible
{\bf Z}-graded $L_0$-submodule of $V_0\otimes {\bf C}[x^{\pm 1}]$.
We know that the $L$-module $V_{{\bf Z}^2}$ has a unique maximal
proper submodule $J_W$ which intersects trivially with $W$. Then we
have the irreducible ${\bf Z}^2$ graded $L$-module
$$
V_{{\bf Z}^2}/J_W.
$$

Now by Theorem 4.3, we have the following result.\vskip2mm

{\em {\bf Theorem 5.1 } $(1)$ If $m_{21}$ is an even integer then
$\OVERLINE{M}_{{\bf Z}^2}^+(\psi,{\bf m}_1,{\bf m}_2)/J_W$ is a
quasi-finite irreducible ${\bf Z}^2$-graded $L$-module for any
exp-polynomial function $\psi$ over $L_0$ and any irreducible ${\bf
Z}$-graded $L_0$-submodule $W$ of $V_0\otimes{\bf C}[x^{\pm 1}]$.

$(2)$ If $m_{21}$ is an odd integer then $ \OVERLINE{M}_{{\bf
Z}^2}^+(\ula ,\psi,{\bf m}_1,{\bf m}_2)/J_W$ is a quasi-finite
irreducible ${\bf Z}^2$-graded $L$-module for any exp-polynomial
function $\psi$ over ${\cal A}$, any finite sequence of nonzero
distinct numbers $\ula =(a_1,\cdots,a_{\nu})$, any finite
dimensional irreducible $sl_2$-modules $V_1,\cdots,V_{\nu}$ and
irreducible {\bf Z}-graded $L_0$-submodule $W$ of $V(\ula
,\psi)\otimes{\bf C}[x^{\pm 1}]$.}

{\bf Remark 5.2 } Since $V_0\otimes{\bf C}[x^{\pm 1}]$ and $V(\ula
,\psi)\otimes{\bf C}[x^{\pm 1}]$ are in general not irreducible
$L_0$-modules, $\OVERLINE{M}_{{\bf Z}^2}^+(\psi,{\bf m}_1,{\bf
m}_2)$ and $ \OVERLINE{M}_{{\bf Z}^2}^+(\ula ,\psi,{\bf m}_1,{\bf
m}_2)$ are in general not irreducible. For example, if $m_{21}$ is
an even integer then we can define an exp-polynomial function $\psi$
over $L_0$ as follows
$$
\psi(t_0^jt^{i{\bf
m}_2})=\frac{(-1)^i+1}{(1-(-1)^jq^{i\alpha})q^{\frac{1}{2}i^2m_{21}m_{22}}},
\quad \psi(\beta)=2,\quad \psi(t_0^{\bar 1}t^{\bf 0})=\frac{1}{2}.
$$
One can  check that $W=v_0\otimes{\bf C}[x^{\pm 2}]$ is an
irreducible {\bf Z}-graded $L_0$-submodule of $v_0\otimes {\bf
C}[x^{\pm 1}]$. Thus the ${\bf Z}^2$-graded $L$-module
$\OVERLINE{M}_{{\bf Z}^2}^+(\psi,{\bf m}_1,{\bf m}_2)$ corresponding
to this function $\psi$ is not irreducible. Suppose $m_{21}$ be an
odd integer. Let $V_1$ be the three dimensional irreducible
$sl_2$-module with the highest weight vector $v_2$. Denote
$E_{21}\cdot v_2$ and $(E_{21})^2\cdot v_2$ by $v_0,v_{-2}$
respectively. Then, for $\ula =(1)$, the exp-polynomial function
$\psi=0$ over ${\cal A}$ and the $sl_2$-module $V_1$, one can see
that
$$
W=\la v_2\otimes x^{2k}\mid k\in{\bf Z}\ra \oplus\la v_{-2}\otimes
x^{2k}\mid k\in{\bf Z}\ra \oplus\la v_0\otimes x^{2k+1}\mid k\in{\bf
Z}\ra ,
$$
is an irreducible {\bf Z}-graded $L_0$-submodule of $V(\ula ,\psi)$.
Thus the corresponding ${Z}^2$-graded $L$-module $
\OVERLINE{M}_{{\bf Z}^2}^+(\ula ,\psi,{\bf m}_1,{\bf m}_2)$ is not
an irreducible module. 
\end{document}